\documentclass[11pt,a4paper]{article}
\usepackage[T1]{fontenc}
\usepackage{lmodern}
\usepackage[margin=1in]{geometry}
\usepackage[numbers,sort&compress]{natbib}
\usepackage{amsmath}
\allowdisplaybreaks
\usepackage{algorithm}
\usepackage{algpseudocode}
\usepackage{setspace}
\usepackage{xcolor}
\usepackage{paralist}
\usepackage{amsthm}
\usepackage{mathtools}

\theoremstyle{plain}%
\newtheorem{theorem}{Theorem}%
\newtheorem{proposition}[theorem]{Proposition}%
\newtheorem{lemma}[theorem]{Lemma}%
\newtheorem{corollary}[theorem]{Corollary}%

\theoremstyle{remark}%
\newtheorem{example}{Example}%

\theoremstyle{definition}%
\newtheorem{definition}{Definition}%
\newtheorem{assumption}{Assumption} 

\usepackage{graphicx}
\usepackage[english]{babel}
\usepackage{amsmath,amssymb,amsfonts}
\usepackage{graphicx}

\DeclareMathOperator*{\argmax}{arg\,max}
\DeclareMathOperator*{\argmin}{arg\,min}
\DeclareMathOperator*{\sign}{{sign}}
\DeclareMathOperator*{\diag}{{diag}}
\newcommand{\Real}{\mathbb R}
\newcommand{\bbb}{{\bf b}}
\newcommand{\ba}{{\bf a}}
\newcommand{\bc}{{\bf c}}
\newcommand{\bd}{{\bf d}}
\newcommand{\be}{{\bf e}}
\newcommand{\bh}{{\bf h}}
\newcommand{\bbf}{{\bf f}}
\newcommand{\bs}{{\bf s}}

\newcommand{\bv}{{\bf v}}
\newcommand{\bu}{{\bf u}}
\newcommand{\bw}{{\bf w}}
\newcommand{\bx}{{\bf x}}
\newcommand{\by}{{\bf y}}
\newcommand{\bz}{{\bf z}}
\newcommand{\bG}{{\bf G}}
\newcommand{\bA}{{\bf A}}
\newcommand{\bB}{{\bf B}}
\newcommand{\bH}{{\bf H}}

\newcommand{\bX}{{\bf X}}
\newcommand{\J}{{\mathcal{J}}}
\newcommand{\M}{{\mathcal{M}}}

\newcommand{\U}{{\mathcal{U}}}
\newcommand{\W}{{\mathcal{W}}}
\newcommand{\balpha}{{\boldsymbol{\alpha}}}
\newcommand{\blbeta}{{\boldsymbol{\beta}}}
\newcommand{\blambda}{{\boldsymbol{\lambda}}}
\newcommand{\bmu}{{\boldsymbol{\mu}}}
\newcommand{\bgamma}{{\boldsymbol{\gamma}}}
\newcommand{\bepsilon}{{\boldsymbol{\epsilon}}}
\newcommand{\bxi}{{\boldsymbol{\xi}}}
\newcommand{\bleta}{{\boldsymbol{\eta}}}
\newcommand{\bpi}{{\boldsymbol{\pi}}}
\newcommand{\bPhi}{{\boldsymbol{\Phi}}}
\newcommand{\bGamma}{{\boldsymbol{\Gamma}}}
\newcommand{\bb}[1]{{\bf #1}}
\newcommand{\norm}[1]{\|#1\|}
\DeclareMathOperator*{\sto}{s.t.}

\usepackage[hidelinks]{hyperref}
\usepackage{authblk}

\begin{document}

\title{Partially Adaptive Robust Optimization for Uncertain Fluid Processing Networks}
\author[1,2]{Roi Ben-Gigi}
\author[2]{Evgeny Shindin}
\author[3]{Shimrit Shtern}
\affil[1]{Faculty of Data and Decision Sciences, Technion -- Israel Institute of Technology, Haifa 3200003, Israel}
\affil[2]{IBM Research, Haifa, Israel}
\affil[3]{Faculty of Engineering, Bar-Ilan University, Ramat-Gan 5290002, Israel}

\date{}
\maketitle

\begin{abstract}

This paper studies the problem of controlling fluid approximations of processing networks under uncertainty. Classical robust formulations for such models rely on static (here-and-now) control decisions that must remain feasible for all realizations of uncertainty, often leading to highly conservative solutions. We propose a partially adaptive robust optimization framework that incorporates observed system information into the control policy at the planning point, enabling less conservative and more responsive decisions while retaining worst-case guarantees.

Our approach begins by modeling the fluid network dynamics as a Separated Continuous Linear Program (SCLP) in the server-effort formulation. We introduce continuous-time extensions of the budgeted uncertainty set and develop three families of uncertainty structures appropriate for fluid-flow systems. For each family, we derive the corresponding robust counterpart. Since these robust counterparts include uncertainty in the objective, they cannot be solved directly by the classical SCLP--Simplex method.

To address this challenge, we develop a hybrid algorithm that integrates an extended SCLP framework with the Saddle-Point Frank--Wolfe (SP--FW) method. This hybrid approach combines tractable continuous-time optimization with efficient saddle-point updates, yielding feasible and near-optimal partially adaptive robust solutions. The resulting methodology advances robust fluid control by enabling uncertainty-aware, dynamically informed decision policies while preserving computational tractability.

\end{abstract}

\section{Introduction}

Processing networks constitute a fundamental modeling paradigm for complex operational systems such as manufacturing lines, customer-service platforms, cloud-computing infrastructures, and sensor networks. The basic entities in such systems are \emph{jobs}, which must be routed, stored, and processed as they move through the network. A typical processing network consists of \emph{buffers}, which temporarily store jobs, and \emph{servers}, which process them. Jobs are partitioned into classes, each associated with a designated buffer; after service completion, a job either exits the system or transitions to another class, joining its corresponding buffer.

The operation of such networks incurs holding and processing costs. The central optimization question is how to schedule, route, and allocate effort in a way that minimizes total operational cost. In principle, one could cast this as a stochastic dynamic optimization problem, but solving it would require exploring an exponential number of scheduling policies, rendering this approach computationally infeasible for large-scale networks.

A widely used and significantly more tractable alternative is to approximate the stochastic network using its {fluid}, or continuous-flow, relaxation. Classical works such as \cite{Weiss1990,Pullan1993} demonstrate that replacing discrete, job-level dynamics with deterministic flows often yields high-quality approximations while enabling the use of continuous optimization methods. In particular, the resulting fluid relaxation of processing networks often takes the form of a Separated Continuous Linear Program (SCLP) and can be efficiently solved using the SCLP--Simplex algorithm of Weiss \cite{Weiss2008, Shindin2021}. The SCLP--Simplex algorithm exploits the structure of these problems, and their close resemblance to classical linear programming models, to obtain optimal piecewise-constant control trajectories. A rigorous theoretical justification for this approach is provided in \cite{Nazarathy2009}, where the authors show that solutions to fluid approximations converge asymptotically, under suitable scaling, to the behavior of the underlying discrete network.

In practice, the arrival rates of jobs and the service rates of servers are seldom known precisely, and may deviate from nominal forecasts. While one may attempt to impose probabilistic models on these uncertainties, such distributions are often difficult to estimate reliably. An alternative approach, adopted in the theory of {Robust Optimization (RO)}, is to describe all uncertain parameters through an {uncertainty set}, representing the set of realizations against which the system must remain protected. RO seeks decisions that (i) remain feasible for all realizations in this set, and (ii) optimize the worst-case value of the objective function; see \cite{BenTal2009} for a comprehensive overview.

Several works apply RO to fluid-processing networks. Bertsimas et al.~\cite{Bertsimas2015} study uncertainty in average processing times and arrival rates under the {effective processing-rate model}, where the decision variables directly control the actual processing rate allocated by each server to each job class. Their formulation employs the {budgeted uncertainty set} of \cite{Bertsimas2004}, which restricts the total number of rate parameters that may simultaneously attain their worst-case deviations. 

Ship et al.~\cite{Ship2022} analyze the {server-effort model}, in which the decision variables specify the fraction of each server's capacity devoted to each flow. They examine both budgeted and general polyhedral uncertainty sets for service rates, together with box uncertainty for arrival rates. Although the effective processing-rate and server-effort models are equivalent in deterministic settings, \cite{Ship2022} shows that their corresponding robust solutions can differ substantially when uncertainty is present.

A related generalization is provided in \cite{Ghate2020}, which develops the robust counterpart of a general Continuous Linear Program (CLP) under budgeted uncertainty. Using a duality-based reformulation, Ghate shows that the robust CLP remains a CLP, thereby encompassing the models of \cite{Bertsimas2015,Ship2022}.

The aforementioned robust formulations for fluid networks employ what the RO literature calls \emph{static} or \emph{here-and-now} optimization, in which all decisions must be chosen before uncertainty is realized. Such policies are often conservative because they must remain feasible for all possible scenarios, even those that never occur.

However, processing networks evolve over time, and by each decision epoch, a substantial portion of the uncertainty, e.g., realized arrivals or queue lengths, has already been observed. RO therefore distinguishes between static (here-and-now) and adaptive (wait-and-see) decisions,
where adaptive decisions may depend on the uncertainty observed up to the time of decision; see, for example, \cite{BenTal2009}.
Adaptive (or ``adjustable'') robust policies can be significantly less conservative than static decisions, and are often better aligned with real operational systems, yet they remain largely unexplored in robust fluid models. 
Models using adaptive policies differ from optimization approaches that re-optimize at each decision point by accounting for adaptivity when making initial decisions. Moreover, re-optimization may not always be possible or may be either computationally or operationally expensive, and thus impracticable.

To obtain tractable adaptive policies, the adjustable robust optimization (ARO) literature commonly restricts the controller to a class of \emph{affine decision rules} (ADRs), in which the decision at time \(t\) is modeled as an affine function of the uncertainty observed up to that time; see \cite{BenTal2004} and \cite[Chapter~14]{BenTal2009}. Such affine parameterizations provide a practical means of incorporating observed system data while keeping the resulting optimization problem computationally manageable. Although ADRs only approximate the fully optimal adaptive policy, they offer a rich and expressive structure that captures essential adaptiveness while preserving tractability.

{Extending this framework to a continuous-time setting, the decision at time 
$t$ becomes an affine function of the entire uncertainty path observed up to 
$t$, so that the controller adapts at every instant. Such continuous adaptation is neither computationally tractable nor consistent with practice, where system information is typically available only at discrete time instances. We therefore consider \textit{partially adaptive policies}, in which the controller adapts only at a predetermined set of observation times. This restricted form of adaptiveness retains much of the expressive power of ADRs while yielding a feasible and efficiently computable solution space.}

{\bf Our Contribution.} We propose a framework that integrates observed system information into the robust control of fluid processing networks, yielding a \emph{partially adaptive affine}  robust solution. As in prior work, we employ robust optimization to ensure feasibility and optimality for all realizations within a prescribed uncertainty set, but we leverage observed state trajectories to reduce conservatism. Within this framework, we formulate three continuous-time budgeted uncertainty sets for fluid-flow systems: a buffer-wise budgeted uncertainty set, as in previous robust fluid-network models, including \cite{Bertsimas2015,Ship2022}; a time-interval budgeted uncertainty set; and a combined buffer-and-time budgeted uncertainty set. For each set, we derive the corresponding robust counterpart. Due to the way uncertainty enters the objective through the partially adaptive controller, the resulting formulations cannot be solved directly using the classical SCLP--Simplex method. To address this limitation, we formulate a saddle-point Frank--Wolfe (SP--FW) scheme tailored to this setting. The scheme combines existing online projection-free ideas with the SCLP--Simplex method, which serves as a linear minimization oracle for the controller subproblem. We establish convergence guarantees for the resulting ergodic sequence and analyze an objective approximation that yields tractable uncertainty-side oracle computations while preserving consistency of the resulting controller.

{\bf Paper structure.} Section~\ref{sec:preliminaries} of this paper introduces the fluid-processing model and establishes the mathematical framework used throughout the paper. Section~\ref{sec:partially_adaptive_control} presents our methodology for constructing a partially adaptive affine controller that incorporates observed system states.  %
{Section~\ref{sec:constraints} introduces the continuous-time budgeted uncertainty sets considered in this paper and derives explicit reformulations of the corresponding robust constraints.} Finally, Section~\ref{sec:Algorithm} casts the resulting optimization problem as a bilinear saddle-point problem, justifies the use of SCLP--Simplex as a linear minimization oracle for the controller subproblem, establishes convergence guarantees for the SP--FW approach, and analyzes an objective approximation used to obtain tractable uncertainty-side oracles.

{\bf Notation.} We use bold symbols $\bx$ and $\bX$ to denote vectors and matrices, respectively. The $j$-th component of $\bx$ is denoted by $x_j$, and the $(j,k)$-component of $\bX$ is denoted by $X_{j,k}$. Row $j$ of matrix $\bX$ is denoted by $\bX_j$. We denote by $\be_i$ the $i$-th vector of the standard basis, and by $\be$ the vector of all ones, with dimension clear from context. For any $l<j$, $\sum_{k=j}^{l} a_k\equiv 0$. For two vectors $\bx,\by\in\mathbb R^n$, $\bx\circ\by$ denotes their element-wise product. Inequalities and equalities involving vectors are understood component-wise. Unless stated otherwise, scalar operations applied to vectors, such as $|\cdot|$, $\operatorname{sign}(\cdot)$, $\min\{\cdot,\cdot\}$, and $\max\{\cdot,\cdot\}$, are interpreted component-wise. We use the convention $\operatorname{sign}(t)=t/|t|$ for $t\neq0$ and $\operatorname{sign}(0)=1$. For $n\in\mathbb N$, we define $[n]=\{1,\ldots,n\}$ and $[n]_0=\{0,1,\ldots,n\}$.

\section{Processing networks and their fluid approximations}\label{sec:preliminaries}

\begin{figure}[htp]
    \centering
    \includegraphics[width=10cm]{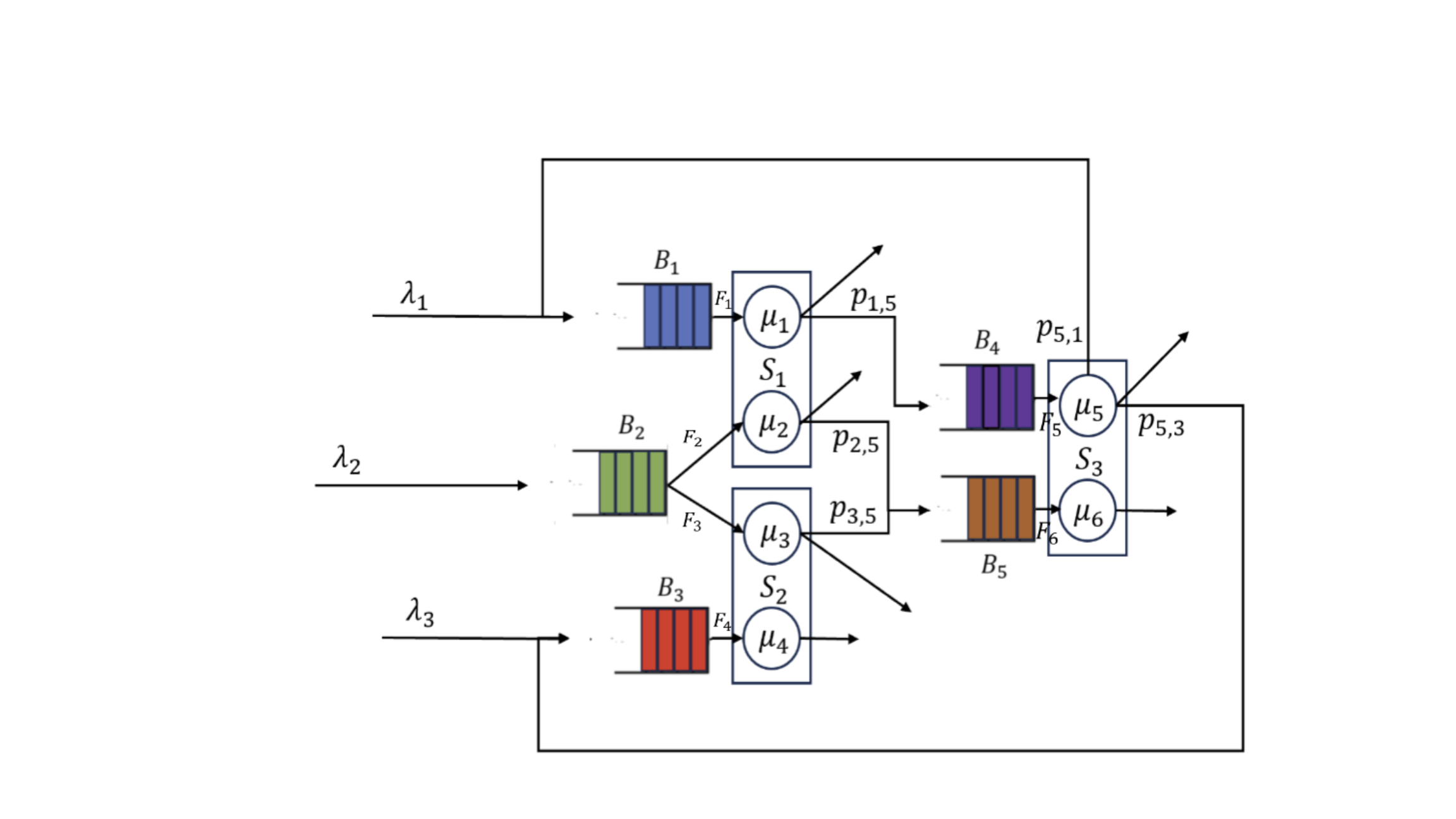}
    \caption{A  processing network with five buffers, three servers, and five job classes. Exogenous arrivals occur at buffers $B_1$, $B_2$, and $B_3$ at rates $\lambda_1,\lambda_2,\lambda_3$, respectively. Rectangles represent jobs, arrows represent possible routing paths, and colors distinguish job classes. After processing, jobs may either leave the system or be reclassified and routed to downstream buffers according to fixed routing probabilities. For instance, jobs in Buffer~2 may be processed by either Server~1 or Server~2 at maximum processing rates $\mu_2$ and $\mu_3$, respectively. Jobs processed by Server~1 leave the system with probability $1-p_{2,5}$ or are routed to Buffer~5 with probability $p_{2,5}$ for further service at Server~3.
    }\label{fig:simple_net}
\end{figure}

Processing networks consist of buffers, which store jobs, servers, which process jobs, and routing rules, which determine how jobs move through the network after processing. The control decisions specify how server capacity is allocated among the available processing activities. Figure~\ref{fig:simple_net} illustrates a simple example with five buffers, three servers, and five job classes. %
To analyze such networks, a common approach is to approximate the stochastic system by a fluid model, %
replacing discrete job counts with continuous fluid levels. In this fluid model each permitted buffer--server pair defines a \emph{flow}, and the control decisions specify the rate at which each flow is processed.

\subsection{Processing-network components}
Consider a fluid network with $K$ buffers, $J$ flows, and $I$ servers. Let $\bx(t) \in \Real^K$ denote the system state at time $t \in [0,T]$, where $x_k(t)$ is the amount of fluid in buffer $k$. Without loss of generality, it is assumed that each flow $j \in [J]$ drains a single buffer $\tilde{k}(j) \in [K]$ and is processed by a single server $\tilde{s}(j) \in [I]$. Each server may process multiple flows.

Let $\blambda\in\Real^K_+$ denote the exogenous arrival rates into the buffers, where $\lambda_k$ is the arrival rate into buffer $k$. Let $\bmu\in\Real^J_+$ denote the maximum processing rates of the flows, where $\mu_j$ is the processing rate attained when server $\tilde s(j)$ dedicates all of its capacity to flow $j$.

Let $\bu(t)\in\Real^J$ denote the vector of processing rates at time $t$, where $u_j(t)$ is the rate at which flow $j$ is processed. For each $j\in[J]$, the processing rate satisfies $0\le u_j(t)\le \mu_j$.
After processing, a fraction $p_{j,k}$ of the processed fluid from flow $j$ is routed to buffer $k\neq \tilde k(j)$, while the remaining fraction $1-\sum_{k\in[K]\setminus\{\tilde k(j)\}}p_{j,k}$ leaves the system.

Two control formulations appear in the literature. In the processing-rate formulation, the control variable is $\bu(t)$. In the server-effort formulation, the control variable is $\bleta(t)$, where $\eta_j(t)$ is the fraction of server $\tilde s(j)$ allocated to flow $j$, so that $\bu(t)=\bmu\circ\bleta(t)$. When $\bmu$ is deterministic, these formulations are equivalent. However, when $\bmu$ is uncertain, they may lead to different decisions. Throughout this work, we adopt the server-effort formulation.

\subsection{Deterministic Fluid Model as an SCLP}

Under the server-effort formulation, the deterministic fluid model admits the
following continuous-time linear optimization formulation:
\begin{subequations}\label{prob:server_effort_model}
\begin{align} 
    \min_{\bx(\cdot), \bleta(\cdot)} \quad & 
        \J(\bleta)
        =
        \int_0^T
        \left[
            \bc^\top \bx(t)
            +
            \bd^\top(\bmu\circ\bleta(t))
        \right] dt  
        \label{eqn:eff0g}\\  
    \text{s.t.} \quad &
        x_k(t)
        =
        \alpha_k + \lambda_k t 
        -
        \int_0^t
        \sum_{j=1}^J
        G_{k,j}\eta_j(s)\,ds,
        &\forall k \in [K],~ t \in [0,T],
        \label{eqn:eff1g}\\
    &
        \sum_{j:\tilde{s}(j)=r}\eta_j(t)\le 1,
        &\forall r\in[I],~t\in[0,T],
        \label{eqn:eff2g}\\
    &
        \bleta(t),\,\bx(t)\ge 0,
        &t\in[0,T].
        \label{eqn:eff3g}
\end{align}
\end{subequations}
Here $\bc\in\Real^K$ and $\bd\in\Real^J$ denote the per-unit holding and processing costs, respectively, and $\balpha=\bx(0)$ is the initial buffer state. We write the objective as $\J(\bleta)$ because the state trajectory $\bx(\cdot)$ is determined by $\bleta(\cdot)$ through the flow-conservation equations.
The matrix $\bG\in\Real^{K\times J}$ is the buffer-flow input-output matrix, where $G_{k,j}$ encodes the effect of serving flow $j$ on the fluid level in buffer $k$ and is defined as
\[
    G_{k,j}
    =
    \begin{cases}
        \mu_j, & k=\tilde{k}(j),\\
        -\mu_j p_{j,k}, & k\ne \tilde{k}(j).
    \end{cases}
\]
Constraint~\eqref{eqn:eff1g} enforces flow conservation, constraint~\eqref{eqn:eff2g} enforces the capacity limit of each server, and constraint~\eqref{eqn:eff3g} ensures non-negativity of the state and control trajectories.

The structure of problem~\eqref{prob:server_effort_model} reveals that the
fluid model is a special case of a \emph{Separated Continuous Linear Program
(SCLP)}, a class introduced by Anderson et al.~\cite{AndersonNashPerold1983}
and subsequently studied in~\cite{Shapiro2001}. In its state-constraint form, an SCLP is an optimization
problem of the form
\begin{equation*}
\begin{aligned}
\min_{\bx(\cdot), {\bu}(\cdot)} \quad
    &
    \int_0^T
    \bc(t)^\top \bx(t)
    +
    \bd(t)^\top \bu(t)
    \,dt
    \\[4pt]
\text{s.t.}\quad
    &
    \int_0^t \bb{B}\bu(s)\,ds
    +\bb{R}\bx(t)=\bb{a}(t), &t\in[0,T]\\[4pt]
    &
    \bb{H}\bu(t)\le \bb{h},
    \qquad
    \bx(t)\ge 0,\quad \bu(t)\ge 0, & t\in[0,T].
\end{aligned}
\end{equation*}
Here $\bx(\cdot)$ denotes the state trajectory and $\bu(\cdot)$ denotes the
control trajectory, and $\bb{a}(t)=\bb{R}\bx(0)+\bar{\bb{a}} t$. In the server-effort model, the buffer levels $\bx(t)$ play
the role of the SCLP state variables, while the service-effort variables
$\bleta(t)$ play the role of the SCLP control variables.

When $\bc(t)$ is constant and $\bd(t)$ is linear, the problem can be solved efficiently by
the SCLP--Simplex algorithm of  \cite{Weiss2008}. In this setting, the
optimal state trajectories are piecewise affine, while the optimal control
trajectories are piecewise constant.

The robust formulation considered in this paper leads to a slightly richer continuous-time linear program. In particular, two additional structural features arise. First, some control variables become active only after prescribed times, leading to constraints of the form
\begin{equation}\label{eq:special_activation_cons}
    \bu_i(t)=0,
    \qquad
    \forall t\in[0,T_i],\quad \forall i\in[n],
\end{equation}
where $\bu_i$ is some subvector of $\bu$ and $T_i$ is a predefined time. 
This can be viewed as replacing $u_i(t)$ by $\chi_i(t)\hat u_i(t)$, where $\chi_i(t)$ is an activation function that equals $0$ on $[0,T_i]$ and $1$ thereafter. Embedding the activation functions into the model coefficients results in matrices $\bB$ and $\bH$ with piecewise constant columns and an objective coefficient vector $\bd(t)$ with piecewise affine components. Alternatively, one may keep $\bB$ and $\bH$ constant and retain $\bd(t)$ as an affine function by imposing the activation constraints explicitly, as in the extended formulation below. Second, the right-hand side $\ba(t)$ of the state-balance constraints is allowed to be piecewise affine. A convenient abstract form for this
extended model is
\begin{equation}\label{eq:extended_sclp_form}
\begin{aligned}[t]
\min_{\bx(\cdot),\bu(\cdot)} \quad
    &
    \int_0^T
    \bc^\top \bx(t)
    +
    \bd(t)^\top \bu(t)
    \,dt
    \\[4pt]
\text{s.t.}\quad
    &
    \int_0^t \bb{B}\bu(s)\,ds+\bb{R}\bx(t)= \ba(t),& t\in [0,T],
    \\[4pt]
    &
    \bb{H}\bu(t)\le \bh%
    ,& t\in [0,T],
    \\[4pt]
    &
    \bu_i(t)=0,
    &
     t\in[0,T_i],\; i\in[n],
    \\[4pt]
    &
    \bx(t)\ge 0,\qquad \bu(t)\ge 0 & t\in [0,T].
\end{aligned}
\end{equation}

We refer to the extended formulation~\eqref{eq:extended_sclp_form} as an SCLP throughout the remainder of the paper. This formulation preserves the essential local structure of a classical SCLP.
Indeed, after collecting all breakpoints arising from the activation times
$T_i$ and from the piecewise affine right-hand sides, the horizon $[0,T]$ can be
partitioned into finitely many intervals. On each such interval, the right-hand
sides and objective coefficients are affine functions of time, while the constraint matrices remain
constant. Thus, between two consecutive breakpoints, the problem has the same
local structure as a classical SCLP, and its optimal control trajectory is expected to be piecewise constant.

The classical SCLP--Simplex algorithm does not directly accommodate activation constraints of the form~\eqref{eq:special_activation_cons}. 
Nevertheless, the preceding observation suggests that the SCLP--Simplex framework can be extended to accommodate such constraints. For the remainder of the paper, we therefore assume the availability of such an extension. Developing and analyzing this extension is outside the scope of this paper and is left for future work.
In Section~\ref{sec:constraints}, we show that the robust constraints arising from the partially adaptive model can be reformulated in the SCLP form introduced above.
We then rely on the assumed extension in Section~\ref{sec:Algorithm}, where the extended SCLP--Simplex framework is used as a linear minimization oracle for the controller subproblem.

\section{Partially Adaptive Controller}\label{sec:partially_adaptive_control}

Consider the server-effort model \eqref{prob:server_effort_model} under uncertainty in the arrival rates $\blambda(t) \in \Real^K_+$. 
We assume that the arrival-rate vector at time $t$ is given by
\begin{equation}\label{eq:lambda_gen}
\blambda(t)=\hat{\blambda}+\tilde{\blambda}\circ \bepsilon(t),
\end{equation}
where $\hat{\blambda},\tilde{\blambda}\in\Real^K_+$ are given vectors representing the nominal arrival rates and the maximum deviations, respectively, and $\bepsilon:[0,T]\to\Real^K$ is an uncertain function belonging to a prescribed uncertainty set $\U$.
We assume that $\U$ is chosen such that
$\blambda(t)\ge 0$
for all $t\in[0,T]$ and all $\bepsilon(\cdot)\in\U$.

Ideally, in this setting, the control at time $t$ would depend on the uncertainty realized up to that time, that is, 
$\bleta(t)\equiv \bleta(t;(\bepsilon(s))_{s\in[0,t]})$. %
{In many practical settings, however, arrival rates are not observable quantities: they are latent parameters of the arrival process, and they influence the system only through their cumulative effect on the state. In such cases, the controller observes the induced fluid levels $\bx(t)$  rather than $\bepsilon(t)$.}
We therefore assume that the control at time $t$ depends on the observed state history, that is, $\bleta(t)\equiv \bleta(t;(\bx(s))_{s\in[0,t]})$. 

This presents two challenges. The first challenge is the computational complexity of finding an optimal or near-optimal controller $\bleta(\cdot)$ within the unrestricted space of general functions. The second challenge stems from the dependency of the controller at each point in time on system states at all preceding points in time. Each of these challenges would make the model intractable (see \cite{BenTal2004} for explanation of intractability in the discrete time setting).

To overcome the first challenge, we restrict the controller $\bleta(\cdot)$ to be an affine function of the state. This is in the spirit of \cite{BenTal2004}, which proposed a framework that uses affine decision rules to adapt to revealed uncertainty in discrete time optimization. 

In many real-world systems, the system state cannot be continuously observed and can only be viewed in predetermined points in time.
Accordingly, for \(p\in\mathbb{N}\), we define a \(p\)-ordered partition of the interval \([0,T]\) as a \(p\)-tuple \((T_i)_{i=0}^{p}\), where
\(0=T_0<T_1<T_2<\cdots<T_p=T\).
We then allow the control at time \(t\) to depend only on the state observations obtained at the preceding observation times \(T_i<t\). This novel feature enables the incorporation of a partially adaptive controller for server resources, accommodating system states revealed over time, while overcoming the second challenge.
To the best of our knowledge, we are the first to propose and integrate such an adaptive approach within robust processing network fluid models.

Before formally introducing the 
{partially adaptive affine controller}, in which control decisions are allowed to depend affinely on
system information observed at a finite number of predetermined time points,
we present a simple concrete example.
The example serves two purposes: first, to illustrate the structural difference between
static (here-and-now) robust control and partially adaptive control; and second, to
demonstrate that even a very limited form of adaptivity—based solely on the observed buffer
state at a single intermediate time—can strictly outperform the conventional static robust
approach in terms of worst-case performance. 
\begin{example}[Motivating Example: Static versus Partially Adaptive Control]
Consider a fluid network with a single buffer and a single server. Under our notation, \(K=J=I=1\), and we omit the indices for convenience. External arrivals enter the buffer at the uncertain rate
\[
\lambda(t)=\hat{\lambda}+\tilde{\lambda}\epsilon(t)
=3+\epsilon(t),
\qquad \epsilon(t)\in[-1,1].
\]
The initial buffer level is \(x(0)=2\), the maximal processing rate is \(\mu=5\), and the processing rate is \(u(t)=5\eta(t)\), where \(\eta(t)\in[0,1]\). The horizon is divided by the observation times \(T_0=0\), \(T_1=1\), and \(T_2=T=3\), and the objective is
\[
J(\eta)=\int_0^3 x(t)\,dt.
\]
The minimum arrival rate \(\lambda(t)=2\) governs robust feasibility, whereas the maximum arrival rate \(\lambda(t)=4\) maximizes the holding cost.

We consider three types of policies: a static controller, which does not change based on the observed state, a general partially adaptive controller which can freely change given an observed state, and a partially adaptive affine controller which can be changed as an affine function of the observed states.
Because processing incurs no cost and earlier processing weakly reduces the buffer level, all three policies considered below have the same optimal controller and state on \([0,1]\):
\[
\eta(t)=
\begin{cases}
1, & t\in[0,2/3],\\[1mm]
2/5, & t\in(2/3,1],
\end{cases}\quad x(t)=
\begin{cases}
2+\int_{0}^t\lambda(s)\,ds-5t, & t\in[0,2/3],\\[1mm]
\int_{0}^t\lambda(s)\,ds-2t, & t\in(2/3,1].
\end{cases}
\]
Under the minimum arrival rate, this controller empties the buffer at \(t=2/3\) and maintains \(x(t)=0\) thereafter. Moreover, %
the cost of this phase is given by $$J_1(\eta)=\int_{0}^1 x(t)\,dt=\int_{0}^1(1-t)\lambda(t)\,dt -1/3.$$

We now move to phase-2.
A static controller cannot depend on the observed value \(x(1)\). Since $x(1)=0$ is possible, in order to protect against $x(t)<0$ for all $t\geq 0$, robust feasibility requires its phase-2 processing rate not to exceed the minimum arrival rate. The optimal static controller is therefore
$\eta^{\mathrm{sta}}(t)=0.4$ %
for all $t\in[1,3]$. Thus, 
Phase-2 cost of the static controller is given by
$$J_2(\eta^{\mathrm{sta}})=\int\limits_1^3x^{\mathrm{sta}}(t)\,dt=\int\limits_1^3 \left(x(1)+\int_1^{t}(\lambda(s)-2)\,ds\right)dt=\int\limits_1^3\int\limits_0^{t}(\lambda(s)-2)\,ds\,dt.$$

A general partially adaptive controller may use the observed value \(x(1)\) to process at the maximal rate until the earliest time at which the buffer could become empty under the minimum arrival rate, and then process at rate \(2\):
\[
\eta^{\mathrm{ada}}(t)=
\begin{cases}
1, & t\in\left(1,1+\dfrac{x(1)}{3}\right],\\[2mm]
\dfrac{2}{5},
& t\in\left(1+\dfrac{x(1)}{3},3\right].
\end{cases}
\]
This controller is robustly feasible, but it is not affine in \(x(1)\), since its switching time depends on the observed state. Denoting $\tau=1+\frac{x(1)}{3}$, the worst case phase-2 cost 
conditional on the observed value \(x(1)\) is given by
\[J_2(\eta^{\mathrm{ada}})=
\int\limits_1^3x^{\mathrm{ada}}(t)\,dt
=\int\limits_1^3(3-s)\lambda(s)\,ds -4+\frac{x(1)^2}{6}.\]

Finally, an optimal partially adaptive affine controller is given by
\[
\eta^{\mathrm{aff}}(t)=
\begin{cases}
\displaystyle \frac{2}{5}+\frac{3}{10}x(1),
& t\in\left(1,\frac{5}{3}\right],\\[2mm]
\displaystyle \frac{2}{5},
& t\in\left(\frac{5}{3},3\right].
\end{cases}
\]
Since \(x(1)\in[0,2]\), this controller satisfies
\(\eta^{\mathrm{aff}}(t)\in[2/5,1]\). Under the minimum future arrival rate, for \(t\in[1,5/3]\),
\[
x(t)
=
x(1)\left(1-\frac{3}{2}(t-1)\right)\geq0,
\]
and \(x(5/3)=0\). Thereafter, the processing rate equals \(2\), so the controller is robustly feasible.

Conditional on the observed value \(x(1)\), the phase-2 cost is then given by
\[
J_2(\eta^{\mathrm{aff}})=\int_1^3x^{\mathrm{aff}}(t)\,dt
=
\int\limits_1^3(3-s)\lambda(s)\,ds -4+\frac{x(1)}{3}.
\]
Plugging in the formula for $x(1)$, adding the phase-1 cost and maximizing over $\lambda(t)\in[2,4]$, the worst-case (maximal) cost $J(\eta)$ of all controllers is attained at $\lambda(\cdot)\equiv 4$  with values
\[
J(\eta^{\mathrm{ada}})
=
J(\eta^{\mathrm{aff}})
=
\frac{19}{3}
<
\frac{29}{3}
=
J(\eta^{\mathrm{sta}}).
\]
Thus, incorporating the state observed at \(T_1\) strictly improves the worst-case objective value. In this example, the partially adaptive affine controller attains the same worst-case value as the more general adaptive controller, without requiring re-optimization at \(T_1\).
\end{example}

This leads us to the following general definition of a partially adaptive affine controller, which we use hereafter.
\begin{definition} \label{def:linear_controller}
    Let $p\in\mathbb{N}$ and let $(T_i)_{i\in[p]_0}$ be a $p$-ordered partition of $[0,T]$.
    We say that $\bleta:[0,T]\rightarrow \Real^{J}$ is a  partially adaptive affine (PAdA) controller with respect to $(T_i)_{i=0}^p$, if %
    \begin{equation} \label{eqn:eta_linear_representation}
    \bleta(t)=\hat{\bleta}^{0}(t)+\sum_{i=1}^{p-1}\bPhi^{i}(t) \bx(T_i),\quad \forall t\in[0,T],
    \end{equation}
where $\hat{\bleta}^{0}: [0,T] \rightarrow \Real^{J}$ is a measurable vector function, and $\bPhi^{i}:[0,T] \rightarrow \Real^{J \times K}$ are measurable matrix functions satisfying $\bPhi^{i}(t) \equiv 0$ for all $t \in [0,T_i]$ and $i\in [p-1]$. 
\end{definition}

\subsection{Affine Dependence on Uncertainty Parameters}

Representation~\eqref{eqn:eta_linear_representation} specifies the controller as an
affine function of the system state evaluated at the observation points. However,
in order to derive %
a tractable robust counterpart, it is necessary to express the
controller explicitly as a function of the uncertainty in the exogenous arrival
rates. Our next goal is therefore to show that the above state-based representation
is equivalent to an affine dependence on the uncertainty process
$\bepsilon(\cdot)$. 

Making this dependence explicit will allow us, in subsequent sections, to analyze
how uncertainty propagates through the system and how it affects the optimization
problem. In particular, we show that {a PAdA controller admits an equivalent} representation
\begin{equation}\label{eq:eta_new0}
\bleta(t)
    = \tilde{\bleta}^{0}(t)
      +\sum_{i=1}^{p-1}\tilde{\bPhi}^{i}(t)
        \int_{T_{i-1}}^{T_i}\bepsilon(s)\,ds,
    \qquad \forall\, t\in[0,T].
\end{equation}

As a preliminary step, we first observe that, without loss of generality, any partially adaptive affine controller of the form~\eqref{eqn:eta_linear_representation} can be expressed as an affine function of the centered states  $\{\tilde{\bx}(T_i)\}_{i\in[p-1]}$, defined as
\begin{align}
\label{eq:def_tildex0}
\tilde{\bx}(T_i)
&\triangleq
\bx(T_i)-\left(\bb{\alpha}+\hat{\blambda} T_i\right) =
-\bG \int_{0}^{T_{i}}\bleta(s)\, ds
+ \int_0^{T_i} \tilde{\blambda} \circ \bepsilon(s)\, ds .
\end{align}
Thus, $\tilde{\bx}(T_i)$ represents the deviation of the observed state from the reference trajectory $\balpha+\hat{\blambda}T_i$.

Using this definition, the linear controller representation \eqref{eqn:eta_linear_representation} can be equivalently rewritten as
\begin{align}
\bleta(t)
&=
\hat{\bleta}^{0}(t)
+ \sum_{i=1}^{p-1} \bPhi^{i}(t)
\bigl(\bb{\alpha} + \hat{\blambda} T_i + \tilde{\bx}(T_i)\bigr) =
\bleta^0(t)
+
\sum_{i=1}^{p-1} \bPhi^{i}(t)\, \tilde{\bx}(T_i),
\label{eta_linear_uncertainty}
\end{align}
where the function $\bleta^{0}:[0,T]\to\Real^J$ is defined as
\begin{equation}
\label{eqn:def_eta_0}
\bleta^0(t)
\triangleq
\hat{\bleta}^{0}(t)
+ \sum_{i=1}^{p-1} \bPhi^{i}(t)
\bigl(\bb{\alpha} + \hat{\blambda} T_i\bigr),
\qquad t\in[0,T].
\end{equation}
Importantly, this reformulation does not impose any additional restrictions on the controller class: for any choice of $\{\bPhi^i(\cdot)\}_{i\in[p-1]}$, the function $\bleta^0(\cdot)$ can be selected independently by an appropriate choice of $\hat{\bleta}^{0}(\cdot)$. %

We now prove the equivalence between \eqref{eta_linear_uncertainty} and \eqref{eq:eta_new0}, which expresses the linear controller $\bleta$ directly in terms of the system’s uncertainty parameters.

\begin{proposition} \label{lem:rep_eta}
\begin{enumerate}[(i)]
\item Let  controller $\bleta$ be of the form \eqref{eta_linear_uncertainty}, \textit{i.e.}, a linear function with given parameter functions $\bleta^{0}$ and $\{{\bPhi}^{i}\}_{i\in[p-1]}$. Then, there exist measurable functions $\tilde{\bleta}^{0}:[0,T]\rightarrow \Real^J$ and $\tilde{\bPhi}^{i}:[0,T]\rightarrow \Real^{J\times K}$ for all $i\in[p-1]$ such that the following equality holds
\begin{align}\label{eq:eta_new}
\bleta(t) &= \tilde{\bleta}^{0}(t)+\sum_{i=1}^{p-1} \tilde{\bPhi}^{i}(t) \int_{T_{i-1}}^{T_{i}} \bepsilon(s) ds, &&\quad \forall t \in [0,T],\\
\tilde{\bPhi}^{i}(t) &\equiv 0, &&\quad \forall i\in [p-1], t \in [0,T_i].\label{eq:zero_phi}
\end{align}
\item  Let $\bleta$ be of the form \eqref{eq:eta_new}, with parameters $\tilde{\bleta}^{0}$ and
$\{\tilde{\bb{\Phi}}^{i}\}_{i=1}^{p-1}$ satisfying \eqref{eq:zero_phi} and
$\tilde{\bb{\Phi}}^{i}_{\cdot,k}(t)=0$ whenever
$\tilde{\lambda}_k=0$.
Consider the functions   
$\bleta^{0}$ and $\{\bb{\Phi}^{i}\}_{i=1}^{p-1}$ such that
for any $t \in [0,T_1]$
\[
\bleta^{0}(t) =\tilde{\bleta}^{0}(t) ,\quad {\bPhi}^{i}(t)   = {\tilde{\bPhi}}^{i}(t) =0, i\in[p-1],
\]
and for $n\in[p-1]$ and %
 $t\in (T_{n},T_{n+1}]$ %
they are constructed recursively according to the following equalities:
\begin{align*}
\bPhi^{i}(t) &  = \tilde{\bPhi}^{i}(t) \diag(\tilde{\blambda})^{\dag} +\bB^{n,i}(t) {\diag}(\tilde{\blambda})^{\dag},&  i= n,n-1,\ldots,1\\
\bleta^{0}(t)&=\tilde{\bleta}^{0}(t) + \bA^{n}(t), %
\end{align*}
where %
\begin{align*}\bB^{n,i}(t)&= \sum_{\ell=i+1}^{n}
\bPhi^{\ell}(t) \bG\int_{0}^{T_{\ell}}  \tilde{\bPhi}^{i}(s) ds -  \sum_{\ell=i+1}^{n} \bPhi^{\ell}(t) 
{\diag}(\tilde{\blambda}),\\
\bA^{n}(t)&=  \sum_{i=1}^{n} \bPhi^{i}(t) \bG \int_{0}^{T_{i}} \tilde{\bleta}^{0}(s) ds.
\end{align*}
Then, \eqref{eta_linear_uncertainty} holds and $\bleta$ is a PAdA controller.
\end{enumerate}
\end{proposition}

\begin{proof}
\begin{enumerate}[(i)]
\item We prove the representation by induction over the intervals
$\{[0,T_i]\}_{i=1}^{p}$, or equivalently over $T=T_i$ for $i\in[p]$.

{\bf Base case:} For $t\in[0,T_1]$ we have $\bPhi^{i}(t)=0$ for all $i\in[p-1]$, since
$T_1\le T_i$ and $\bPhi^{i}$ vanishes on $[0,T_i]$. Hence \eqref{eta_linear_uncertainty}
gives $\bleta(t)=\bleta^{0}(t)$ on $[0,T_1]$, and defining
$\tilde{\bleta}^0(t)\equiv\bleta^{0}(t)$ and $\tilde{\bPhi}^i(t)=0$ for all $t\in[0,T_1]$
and $i\in[p-1]$ satisfies \eqref{eq:eta_new} and \eqref{eq:zero_phi}.

{\bf Induction step:} Assume that for some $n\in[p-1]$ there exist functions
$\tilde{\bleta}^{0}(\cdot)$, $\{\tilde{\bPhi}^{i}(\cdot)\}_{i=1}^{p-1}$ defined over
$\cup_{i=1}^n(T_{i-1},T_i]=(0,T_{n}]$ such that \eqref{eq:eta_new} and
\eqref{eq:zero_phi} hold for $T=T_n$. We show that they can be extended to
$(T_{n},T_{n+1}]$ so that \eqref{eq:eta_new} and \eqref{eq:zero_phi} continue to hold
for $T=T_{n+1}$.

We start by setting $\tilde{\bPhi}^{i}(t)=0$ for all $i=n+1,\ldots,p-1$ and
$t\in(T_{n},T_{n+1}]$. Plugging \eqref{eq:def_tildex0} for
$\tilde{\bx}(T_1),\ldots,\tilde{\bx}(T_{p-1})$ into \eqref{eta_linear_uncertainty} and
using $\bPhi^i(t)=0$ for $t\in[0,T_i]$ yields
\begin{align*}
\bleta(t) = \bleta^{0}(t) + \sum_{i=1}^{n}\bPhi^{i}(t)\left(-\bG\int_{0}^{T_{i}}\bleta(s)\,ds + \int_0^{T_i}\tilde{\blambda}\circ\bepsilon(s)\,ds\right),\quad t\in(T_{n},T_{n+1}].
\end{align*}
Applying the induction hypothesis, i.e.\ substituting \eqref{eq:eta_new} for
$t\in[0,T_n]$, and using $\tilde{\bPhi}^{i}(t)=0$ for $t\in[0,T_i]$, $i\in[n]$, which
truncates the inner sum at $\ell=i-1$, we obtain, for all $t\in(T_{n},T_{n+1}]$,
\begin{align}
\bleta(t) = &\ \bleta^{0}(t) + \sum_{i=1}^{n}\bPhi^{i}(t)\int_0^{T_i}\tilde{\blambda}\circ\bepsilon(r)\,dr - \label{eqn:eta_after_sub}\\
&\ \sum_{i=1}^{n}\bPhi^{i}(t)\,\bG\int_{0}^{T_{i}}\tilde{\bleta}^{0}(s)+\sum_{\ell=1}^{i-1}\tilde{\bPhi}^{\ell}(s)\left(\int_{T_{\ell-1}}^{T_{\ell}}\bepsilon(r)\,dr\right)ds. \nonumber
\end{align}
By linearity of sums and integrals,
\begin{align*}
\sum_{i=1}^{n}\bPhi^{i}(t)\,\bG\int_{0}^{T_{i}}\sum_{\ell=1}^{i-1}\tilde{\bPhi}^{\ell}(s)&\left(\int_{T_{\ell-1}}^{T_{\ell}}\bepsilon(\tau)d\tau\right)ds =
\sum_{i=1}^{n-1}\sum_{\ell=i+1}^{n}\bPhi^{\ell}(t)\,\bG\int_{0}^{T_{\ell}}\tilde{\bPhi}^{i}(s)\,ds\smashoperator{\int_{T_{i-1}}^{T_{i}}}\bepsilon(r)\,dr,\\
\sum_{i=1}^{n}\bPhi^{i}(t)\smashoperator{\int_0^{T_i}}\tilde{\blambda}\circ\bepsilon(r)\,dr =&\sum_{i=1}^{n}\bPhi^{i}(t)\diag(\tilde{\blambda})\sum_{\ell=1}^{i}\int_{T_{\ell-1}}^{T_{\ell}}\!\bepsilon(r)\,dr = \sum_{\ell=1}^{n}\sum_{i=\ell}^{n}\bPhi^{i}(t)\diag(\tilde{\blambda})\smashoperator{\int_{T_{\ell-1}}^{T_{\ell}}}\!\bepsilon(r)\,dr,
\end{align*}
so that, by the definitions of $\bA^{n}$ and $\bB^{n,\ell}$, \eqref{eqn:eta_after_sub}
is equivalent to
\begin{align}\label{eq:recursion}
\hspace{-8pt}\bleta(t) = \bleta^{0}(t) - \bA^{n}(t) + \sum_{\ell=1}^{n}\hspace{-2pt}\left(-\bb{B}^{n,\ell}(t)+\bPhi^{\ell}(t)\diag(\tilde{\blambda})\right)\hspace{-2pt}\smashoperator{\int_{T_{\ell-1}}^{T_\ell}}\!\bepsilon(r)\,dr,\,\,\forall t\in(T_{n},T_{n+1}].
\end{align}
Therefore, extending $\{\tilde{\bPhi}^i\}_{i=1}^{p-1}$ and $\tilde{\bleta}^0$ to
$t\in(T_{n},T_{n+1}]$ by
\begin{align}\label{phi_tilde}
\tilde{\bPhi}^{i}(t) &= \begin{cases} -\bb{B}^{n,i}(t)+\bPhi^{i}(t)\diag(\tilde{\blambda}), & i\in[n],\\ \mathbf{0}, & i=n+1,\ldots,p-1,\end{cases}\\
\tilde{\bleta}^{0}(t) &= \bleta^{0}(t) - \bA^{n}(t),\label{eq:eta_tilde}
\end{align}
yields the desired result and completes the induction.

\item Let now $\{\tilde{\bPhi}^i\}_{i=1}^{p-1}$ and $\tilde{\bleta}^0$ be
given, with $\tilde{\bPhi}^{i}_{\cdot,k}(t)=0$ whenever $\tilde{\lambda}_k=0$. We
construct corresponding $\{\bPhi^{i}\}_{i=1}^{p-1}$ and $\bleta^0$ by inverting
\eqref{phi_tilde} and \eqref{eq:eta_tilde} to obtain \eqref{eq:recursion}.

For $t\in[0,T_1]$, the base case gives $\bleta^0(t)=\tilde{\bleta}^0(t)$ and
$\bPhi^i(t)=0$ for all $i\in[p-1]$. Fix $n\in[p-1]$ and consider
$t\in(T_{n},T_{n+1}]$. 
\[
\bPhi^{i}(t) = \left(\tilde{\bPhi}^{i}(t)+\bb{B}^{n,i}(t)\right)\diag(\tilde{\blambda})^{\dag}.
\]
By its definition, $\bB^{n,i}(t)$ depends only on $\{\bPhi^{\ell}(t)\}_{\ell=i+1}^{n}$
on $(T_{n},T_{n+1}]$ and on the given $\tilde{\bPhi}^{i}$. Hence, proceeding in
descending order $i=n,n-1,\ldots,1$, the above relation determines
$\{\bPhi^{i}(t)\}_{i=1}^{n}$ on $(T_{n},T_{n+1}]$; the remaining functions
$\{\bPhi^{i}(t)\}_{i=n+1}^{p-1}$ are set to $\mathbf{0}$ there. Using
$\{\bPhi^{i}(t)\}_{i=1}^{{n}}$ and $\tilde{\bleta}^0$ we then evaluate $\bA^{n}(t)$ and
obtain $\bleta^{0}(t)=\tilde{\bleta}^{0}(t)+\bA^{n}(t)$ from \eqref{eq:eta_tilde}.
Thus, \eqref{eq:recursion} holds, implying that \eqref{eta_linear_uncertainty} holds as well. Determining $\hat{\bleta}^0$ via \eqref{eqn:def_eta_0}, we conclude that $\bleta$ is of the form~\eqref{eqn:eta_linear_representation}.
\end{enumerate}\vspace{-10pt}
\end{proof}

\subsection{Robust Server Effort Model with PAdA Controller}

To allow partial adaptivity of the controller $\bleta$, we substitute the affine representation~\eqref{eq:eta_new} into the optimization model. Moreover, we eliminate the state variable $\bx$ from problem~\eqref{prob:server_effort_model} and reformulate it as an optimization problem in terms of the controller parameters only.

In order to  eliminate state variables $\{x_k(t)\}_{k=1}^K$ from the objective functional, we can use equation \eqref{eqn:eff1g} to express $x_k(t)$ in terms of $\{\eta_j(t)\}$, and then substitute the result in \eqref{eqn:eff0g}. Using the definition of $\blambda(t)$ at \eqref{eq:lambda_gen} and omitting constant terms involving $\alpha_k$ and $\hat{\lambda}_k$ we obtain %
\begin{equation}\label{eqn:obj_term_raw_bas}
F(\bleta,\bb\epsilon)= \int_0^T \left[-\int_0^t \mathbf{c}^\top \bG   \bleta(s) ds  +  \int_{0}^{t} [\tilde{\blambda} \circ\bc]^\top \bepsilon(s) ds +  \mathbf{d}^\top(\bmu \circ \bleta(t))\right] dt.
\end{equation}

Plugging \eqref{eq:eta_new} into \eqref{eqn:obj_term_raw_bas} and applying Fubini's theorem to the terms involving $\bleta$
yields the objective function
\begin{align}
 \hspace{-10pt} F(\tilde{\bleta}^{0}, \{\tilde{\bPhi}^{i}\}, \bb{\epsilon})
&=\; \smashoperator{\int_0^T} \left[
\left[
-(T-s)\mathbf{G}^\top\mathbf{c}
+(\mathbf{d}\circ\bmu)
\right]^\top
\tilde{\bleta}^{0}(s)\, + \int_0^s
(\tilde{\blambda}\circ\mathbf{c})^\top \bb{\epsilon}(t)\,dt\right]\,ds \label{eqn:obj_term_raw}\\
&+ \sum_{i=1}^{p-1} \int_0^T  
\left[
-(T-s)\mathbf{G}^\top\mathbf{c}
+(\mathbf{d}\circ\bmu)
\right]^\top
\tilde{\bPhi}^{i}(s) 
\,ds \left(\int_{T_{i-1}}^{T_i}\bb{\epsilon}(t)\,dt\right).\nonumber
\end{align}

Treating the state variables $\{x_k(t)\}_{k=1}^K$ in problem \eqref{prob:server_effort_model} as slack variables,  we can replace the constraints of type \eqref{eqn:eff1g} with the equivalent constraints %
\begin{equation} \tag{\ref{eqn:eff1g}'} \label{constraint_x_slack}
\int_0^t  \sum_{j=1}^J G_{k,j} \eta_{j}(s) d s \leq \alpha_k + \int_{0}^t {\lambda}_k (s) ds.
\end{equation}
where the term $t\lambda_k $, which considers a constant known $\blambda$,
is replaced by the integral of the unknown $\lambda_k(s)$. 

Thus, plugging  \eqref{eq:eta_new} and \eqref{eq:zero_phi} into the resulting \eqref{prob:server_effort_model}, the robust counterpart of the server-effort model associated with a partially adaptive affine controller is given by %

\begin{subequations}\label{prob:server_effort_model_only_uncertainty_general}
\begin{alignat}{2}
& \hspace{-3em}\min_{\tilde{\bleta}^{0},\,\{\tilde{\bPhi}^{i}\}_{i\in[p-1]}}
\max_{\bepsilon\in\U}
F(\tilde{\bleta}^{0},\{\tilde{\bPhi}^i\}_{i\in[p-1]},\bepsilon) 
\label{prob:server_effort_model_only_uncertainty_general:a}
\\[0.4em]
\text{s.t.}\qquad
&\hspace{-1em}
\sum_{i=1}^{p-1}
\int_{T_i}^t\!\bG_k\tilde{\bPhi}^{i}(\tau)\,d\tau
\left(\int_{T_{i-1}}^{T_i}\!\!\!\bepsilon(s)\,ds\right)+
\notag\\
&\hspace{-0.9em}
\int_0^t\!\!\bG_k\tilde{\bleta}^{0}(\tau)\,d\tau
\le
\alpha_k+\hat{\lambda}_k t
+\!\!\int_0^t\!\tilde{\lambda}_k\epsilon_k(s)\,ds,
&&\hspace{-0.7em}\forall k\in[K],\;t\in[0,T],\;\bepsilon\in\U,
\label{prob:server_effort_model_only_uncertainty_general:b}
\\[0.4em]
&\hspace{-1em}
\sum_{j:\,\tilde{s}(j)=r}
\left[
\tilde{\eta}^{0}_j(t)
+\sum_{i=1}^{p-1}
\tilde{\bPhi}^{i}_j(t)
\int_{T_{i-1}}^{T_i}\!\!\!\bepsilon(s)\,ds
\right]
\le 1,
\quad
&&\forall r\in[I],\;t\in[0,T],\;\bepsilon\in\U,
\label{prob:server_effort_model_only_uncertainty_general:c}
\\[0.4em]
&\hspace{-1em}
\tilde{\bleta}^{0}(t)
+\sum_{i=1}^{p-1}
\tilde{\bPhi}^{i}(t)
\int_{T_{i-1}}^{T_i}\!\!\!\bepsilon(s)\,ds
\ge\bb{0},
\quad
&&\forall t\in[0,T],\;\bepsilon\in\U,
\label{prob:server_effort_model_only_uncertainty_general:d}
\\[0.4em]
&\hspace{-1em}
\tilde{\bPhi}^{i}(t)=\bb{0},
\quad
&&\forall t\in[0,T_i],\;i\in[p-1].
\label{prob:server_effort_model_only_uncertainty_general:e}
\end{alignat}
\end{subequations}
In the next sections, we discuss solving this model for various uncertainty sets $\U$ through reformulations and first-order methods.

\section{Robust Constraint Reformulations}\label{sec:constraints}
In this section, 
we reformulate the constraints of problem \eqref{prob:server_effort_model_only_uncertainty_general} as SCLP constraints for different types of uncertainty sets $\U$. All the  uncertainty sets we consider here are generalizations of the {\em{budgeted}} uncertainty set \cite{Bertsimas2004}. Generally, a budgeted uncertainty set limits the amount of uncertainty spread between different uncertain parameters according to a fixed budget. Intuitively, this type of uncertainty set is motivated by a reasonable and intuitive expectation that the worst-case realization would not occur in all uncertain parameters simultaneously, and hence the too-conservative solution that protects from such a worst-case scenario is not necessary. Originating in \cite{Bertsimas2004}, budgeted uncertainty sets were introduced in the context of LP problems, and shown to ensure probabilistic guarantees on the solution, under 
the assumption that the uncertain noise variables are independent and symmetrically distributed random variables in [-1,1]. 

The budgeted uncertainty sets we present and analyze in this work differ from one another in how the budget is distributed. First, in subsection \ref{sec:uncert_in_all_buffers}, we use a budgeted uncertainty set that bounds the sum of uncertainty in the exogenous input rates to all buffers at each point in time, similarly to \cite{Bertsimas2015,Ship2022},  while in subsection \ref{sec:uncert_over_time} the budgeted uncertainty limits the total deviation of each of the input rates to a single buffer over a  predefined time interval.
Finally in subsection \ref{sec:uncertainty_set_3_combined} we consider an uncertainty set that generalizes the two aforementioned sets limiting the uncertainty over both buffers and time. 

Note that all constraints of problem \eqref{prob:server_effort_model_only_uncertainty_general} can be characterized as having the following function of the uncertain parameter $\bepsilon$:
\begin{equation}\label{eq:general_cons_epsilon_ext}
\sum_{i=1}^{l-1} {\bf b}^i(t)^\top \int_{T_{i-1}}^{T_i} \bepsilon(s) ds -   \int_{T_{l-1}}^{t} \tilde{\lambda}_{k}^{0}\bepsilon_k(s) ds \leq z(t) \quad t \in (T_{l-1},T_{l}], \; l \in [p],\; \forall \bepsilon\in \U,
\end{equation}
for some $k\in[K]$, where $\tilde{\lambda}_{k}^{0}$ is equal either to $[\tilde{\blambda}]_k$ or to $0$, and $\{b_i(\cdot)\}$ are general functions, which may depend on the decision variables. Since the feasibility of a solution must hold for all $t \in [0,T]$ and for all uncertainty
realizations $\bepsilon \in \U$, we fix an arbitrary $t \in [0,T]$ and study the resulting worst-case value of
the left-hand side of \eqref{eq:general_cons_epsilon_ext} over all admissible uncertainty
realizations.

\subsection{Budgeted Uncertainty Over Buffers} \label{sec:uncert_in_all_buffers}
In this section, we study a budgeted uncertainty set that constrains the aggregate perturbation of the exogenous input rates at any given time.
For any $\mathcal{T}\subseteq [0,T]$, let
\[
\M^{\mathcal{T},K}
\coloneqq
\left\{
\bepsilon:\mathcal{T}\to\mathbb{R}^{K}
\;\middle|\;
\bepsilon \text{ is Lebesgue measurable}
\right\}.
\]
We describe the uncertainty set $\U_1$ as

\begin{equation}\label{def:uncertainty_set_1}
\U_1
    = \left\{
        \bepsilon(\cdot)\in \M^{[0,T],K} \;\middle|\;
        \|\bepsilon(t)\|_{\infty} \le 1,\;
        \|\bepsilon(t)\|_{1} \le \Gamma(t),
        \; \forall\, t \in [0,T]
      \right\}.
\end{equation}

We use the following standing assumption on the budget function $\Gamma(\cdot)$.
\begin{assumption}\label{ass:Gamma}
    The function $\Gamma:[0,T]\rightarrow [0,K]$, %
    is piecewise constant with a finite number of pieces $q$.
\end{assumption}
Indeed, the assumption that the image of $\Gamma$ lies in $[0,K]$ can be made without loss of generality, since negative $\Gamma(t)$ leads to an empty set $\U_1$, and for each $t\in[0,T]$ such that $\Gamma(t)>K$ taking $\Gamma(t)=K$ produces the same set $\U_1$.

To ensure that constraint~\eqref{eq:general_cons_epsilon_ext} holds for all
$\bepsilon \in \U_1$, we must maximize the left-hand side of
\eqref{eq:general_cons_epsilon_ext} over all admissible uncertainty realizations.
Because the uncertainty-set constraints in $\U_1$ are separable across time,
this maximization decomposes into independent subproblems defined over disjoint time
intervals corresponding to the components of the sum in
\eqref{eq:general_cons_epsilon_ext}. Specifically, for a given controller, verifying feasibility for every
$t \in (T_{l-1},T_l]$ for $l \in [p]$, reduces to solving $l$ optimization problems, 
where 
\begin{align}
\label{eqn:seperable_prob_i}
      \max \quad& {\bf b}^i(t)^\top \int_{T_{i-1}}^{T_i} \bepsilon(s) ds  \\ 
 \sto \quad &   \sum_{k=1}^K |[\bepsilon]_k(s)|  \leq \Gamma(s), \quad &&\forall s \in (T_{i-1},T_i], \nonumber \\ 
 \quad &  |[\bepsilon]_k(s)|  \leq 1, \quad && \forall s \in (T_{i-1},T_i], \quad k \in [K], \nonumber 
 \end{align} 
computes the worst case for interval $(T_{i-1},T_i]$ for all $i \in [l-1]$
, and 
\begin{align} \label{eqn:seperable_prob_t}
      \max \quad&  - \int_{T_{l-1}}^t \tilde{\lambda}_{k}^{0} [\bepsilon]_k(s) ds  \\
\sto \quad &   \sum_{k=1}^K |[\bepsilon]_k(s)|  \leq \Gamma(s), \quad &&\forall s \in (T_{l-1},T_{l}],\nonumber \\ 
 \quad &  |[\bepsilon]_k(s)|  \leq 1, \quad  &&\forall s \in (T_{l-1},T_{l}], \quad k\in [K]. \nonumber 
  \end{align}
computes the worst case on the interval
$(T_{l-1},t]$, and ensuring that the sum of the resulting worst-case values remains smaller or equal to $z(t)$. These subproblems intersect only at the single time point $T_i$, which does
not affect the objective value and therefore does not impact the validity of the
decomposition.
 
Note that problem~\eqref{eqn:seperable_prob_t} admits a trivial optimal solution given by 
\[
\epsilon_k(s)=-\min\{1,\Gamma(s)\},\qquad  \epsilon_{k'}(s)=0,\;\forall s\in(T_{l-1},T_{l}],\ k'\in[K]\setminus \{k\},
\]
which yields the objective value
\(\tilde{\lambda}_{k}^{0} \int_{T_{l-1}}^{t} \min\{\Gamma(s),1\}\, ds .
\)

Furthermore, for any relevant index~$i$, problem~\eqref{eqn:seperable_prob_i} is equivalent to the following optimization problem:
\begin{equation}\label{eq:seperable_prob_i_abs}
\begin{aligned}[t]
\max \quad 
& \bigl|{\bf b}^i(t)\bigr|^\top \int_{T_{i-1}}^{T_i} \bepsilon(s)\, ds \\[0.3em]
\sto \quad 
& \be^\top \bepsilon(s) \le \Gamma(s), 
\qquad &&\forall s \in (T_{i-1},T_i], \\
& \bb{0} \le \bepsilon(s) \le \be, 
\qquad && \forall s \in (T_{i-1},T_i].
\end{aligned}
\end{equation}
This equivalence follows from the observation that any feasible solution $\bepsilon(\cdot)$ to problem \eqref{eq:seperable_prob_i_abs} can be transformed into a feasible solution ${\sign}({\bf b}^i(t)) \bepsilon(\cdot)$ to problem \eqref{eqn:seperable_prob_i}, and vice versa, $|\bepsilon(\cdot)|$ is feasible to problem \eqref{eqn:seperable_prob_i}, for any feasible solution $\bepsilon(\cdot)$ to \eqref{eq:seperable_prob_i_abs}. In both cases the objective is given by $|{\bf b}^i(t)|^T \int_{T_{i-1}}^{T_i}|\bepsilon(s)|ds$, hence both problems must have same optimal value as problem \eqref{eq:seperable_prob_i_abs}. Hence, we focus on solving problem \eqref{eq:seperable_prob_i_abs}, %
and show it suffices to solve a finite collection of simpler linear programs with finitely many variables and constraints.
These subproblems can then be solved efficiently using classical robust optimization 
methods for linear programs, and their solutions collectively yield the desired 
solution to \eqref{eqn:seperable_prob_i}.

For every $n \in \{0,1,\dots,K\}$, we define $I_n(s)$ to be an indicator function satisfying 
\begin{align*}I_n(s) = \begin{cases} 1, & \mbox{if } \Gamma(s) \in [n, n+1){,}  \\
 0{,} & \text{otherwise}. \end{cases}
 \end{align*}
Denoting $R_{i,n}= (T_{i-1},T_i] \cap \{s\in[0,T]: I_n(s) = 1\}$, the integral over $R_{i,n}$ is defined as 
\begin{equation}
\int_{R_{i,n}} \bepsilon(s) ds \equiv \int_{T_{i-1}}^{T_i} \bepsilon(s) I_n(s) ds.\label{def:integral_over_R}
\end{equation}
Using this definition, the objective function of \eqref{eq:seperable_prob_i_abs} can be rewritten as 
\[ \bigl|{\bf b}^i(t)\bigr|^\top \int_{T_{i-1}}^{T_i} \bepsilon(s)\, ds = |{\bf b}^i(t)|^\top \sum_{n=0}^{K}\int_{R_{i,n}} \bepsilon(s) ds   \] 
Note that $R_{i,n}, n=0,1,\dots, K$ defines a disjoint partition of the interval $(T_{i-1}, T_i]$ and thus, by the separability of the constraints for each $s$, we can separate \eqref{eq:seperable_prob_i_abs} into $K+1$ independent sub-problems of the form
 \begin{align}
\label{eqn:un_primal_sep}
     \max  \quad &|{\bf b}^i(t)|^\top \int_{R_{i,n}} \bepsilon(s) ds  \\ 
 \sto \quad &   \sum_{k=1}^K \bepsilon_k(s)  \leq \Gamma(s), \quad \forall s \in R_{i,n},\nonumber \\ 
 \quad & 0 \leq \bepsilon_k(s)  \leq 1, \quad  \forall s \in R_{i,n}. \nonumber 
  \end{align} 
Clearly, problems associated with $i,n$ for which $R_{i,n}$ has measure zero can be excluded, as they have a value of zero, and so do not affect the total objective value. Therefore, without loss of generality we assume that $R_{i,n} $ is not of measure zero. We also define the duration of $R_{i,n}$ as $S_{i,n} = \int_{R_{i,n}} 1 ds>0$, as well as the average $\Gamma(s)$ over $R_{i,n}$, given by
\[
\Gamma_{i,n} =
\dfrac{\displaystyle \int_{R_{i,n}} \Gamma(s)\, ds}{S_{i,n}}.
\] 

We now present two lemmas which together establish that the 
infinite-dimensional structure of problem~\eqref{eqn:un_primal_sep}, which includes 
infinitely many constraints and function-valued variables, can be reduced to 
a finite-dimensional linear program, enabling computation of its optimal value.
 Lemma \ref{lem:optimal_sol} identifies a potential solution to problem \eqref{eqn:un_primal_sep}, while Lemma \ref{lem:fixed_budget_problem}
leverages Lemma \ref{lem:optimal_sol} to show that replacing the right-hand side $\Gamma(s)$ of the first constraint of \eqref{eqn:un_primal_sep}, with $\Gamma_{i,n}$, results in  the same
objective function
value. Thus,
problem \eqref{eqn:un_primal_sep} can be
replaced by this new
fixed-budget version,
which we use to derive
a finite reformulation of the constraint \eqref{eq:general_cons_epsilon_ext}.

\begin{definition}\label{def:avop}
For any $\bx \in \mathbf{R}^n$, the \emph{absolute value ordering permutation}
$\bpi(\bx)$ is the permutation of $[n]$ that sorts the components of $\bx$ by
non-increasing absolute value, breaking ties in favour of the smaller index.
Formally, $\bpi(\bx)$ is the unique permutation satisfying
\[
|x_{\pi(\bx)_i}| \;\geq\; |x_{\pi(\bx)_j}|
\qquad \text{and} \qquad
\bigl(|x_{\pi(\bx)_i}| = |x_{\pi(\bx)_j}| \;\Longrightarrow\; \pi(\bx)_i < \pi(\bx)_j\bigr)
\]
for all $i,j\in[n]$ with $i<j$.
\end{definition}

\begin{lemma} \label{lem:optimal_sol}
Let $l\in [p]$, $t\in(T_{l-1},T_{l}]$, $i\in [l-1]$, $\tilde{\bpi}=\bpi({\bf b}^i(t))$ and $n\in [K]_0$. Then,
the vector-valued function $ \tilde{\bepsilon}:R_{i,n}\rightarrow [0,1]^K $, given for every $s \in R_{i,n}$ by
\begin{equation}\label{def:sol_def}
\tilde{\epsilon}_k(s) = \begin{cases} 
1, & \text{if } k \in \{\tilde{\pi}_1, \dots, \tilde{\pi}_n\}, \\ 
\Gamma(s) - n, & \text{if } n<K \text{ and } k = \tilde{\pi}_{n+1}, \\
0, & \text{otherwise},
\end{cases}
\end{equation}
is an optimal solution to the problem \eqref{eqn:un_primal_sep}.\end{lemma}
\begin{proof}
Fix  $\tilde{s} \in R_{i,n}$, and consider the optimization problem 
 \begin{align}
\label{eqn:un_primal_sep_s}
      \max \quad & |{\bf b}^i(t)|^\top \bb{\psi}  \\ 
 \sto \quad &   \sum_{k=1}^K [\bb{\psi}]_k  \leq \Gamma(\tilde{s}),  \nonumber\\ 
 \quad & 0 \leq [\bb{\psi}]_k  \leq 1. \nonumber 
  \end{align}

Following \cite{Bertsimas2004}, the solution of problem \eqref{eqn:un_primal_sep_s} is given by $\psi = \tilde{\bepsilon}(\tilde{s})$, where $\tilde{\bepsilon}(\tilde{s})$ is defined according to \eqref{def:sol_def}.
The optimality of $\tilde{\bepsilon}(s)$ for problem~\eqref{eqn:un_primal_sep_s} for each $s\in R_{i,n}$, implies that $\tilde{\bepsilon}(\cdot)$ is a feasible solution to \eqref{eqn:un_primal_sep}.
 Optimality of $\tilde{\bepsilon}(s)$ for the latter therefore gives
$|{\bf b}^i(t)|^\top \tilde{\bepsilon}(s) \geq |{\bf b}^i(t)|^\top \bepsilon(s), \, \forall s \in R_{i,n}$,
for any feasible $\bepsilon(\cdot)$ to problem \eqref{eqn:un_primal_sep}. The integral operation preserves inequality sign, and thus $ \int_{R_{i,n}} |{\bf b}^i(t)|^\top \tilde{\bepsilon}(s) ds \geq \int_{R_{i,n}} |{\bf b}^i(t)|^\top \bepsilon(s) ds $ for any feasible $\bepsilon(\cdot)$, which proves the optimality of feasible solution $\tilde{\bepsilon}(s)$ to problem \eqref{eqn:un_primal_sep}.
\end{proof}

\begin{lemma} \label{lem:fixed_budget_problem}
Let 
\begin{align*}
\U_{i,n}   = \{\bepsilon(\cdot)\in \M^{R_{i,n},K} %
\, | \, 0 \leq \bepsilon(t) \leq \be,\;  \norm{\bepsilon(t) }_1\leq \Gamma_{i,n},  t\in R_{i,n}\}.
\end{align*}
Then the optimization problem  \begin{align}
\label{eqn:un_primal_sep_fixed_budget_b_abs}
    \max_{\bepsilon\in \U_{i,n}}\quad&    |{\bf b}^i(t)|^\top \int_{R_{i,n}} \bepsilon(s) ds  %
  \end{align} 
 has the same optimal value as problem \eqref{eqn:un_primal_sep}.
\end{lemma}
\begin{proof}
Applying Lemma \ref{lem:optimal_sol} to problem \eqref{eqn:un_primal_sep_fixed_budget_b_abs} with the function $\Gamma(s) \equiv \Gamma_{i,n}$ over $R_{i,n}$, and using the definition of  $\Gamma_{i,n}\in[n, n+1)$, we conclude that there exists an optimal solution to \eqref{eqn:un_primal_sep_fixed_budget_b_abs} which is a constant function (independent of $s$) 
$\hat{\bepsilon}:R_{i,n}\rightarrow \Real^K$, %
given by
\[
\hat{\epsilon}_k(s) = \begin{cases} 1, & k \in \{\tilde{\pi}_1,\dots,\tilde{\pi}_n\}, \\ \Gamma_{i,n}-n, & n<K \text{ and } k = \tilde{\pi}_{n+1} ,
\\ 0, & \mbox{else},  
\end{cases}
\]
where $\tilde{\bpi}=\bpi({\bf b}^i(t))$ is the absolute order permutation of ${\bf b}^i(t)$.
Using Lemma \ref{lem:optimal_sol} we know that $\tilde{\bepsilon}(s)$, defined as in \eqref{def:sol_def}, is an optimal solution to problem $\eqref{eqn:un_primal_sep}$. 
For $n=K$, both optimal solutions are equal to $\be$ on $R_{i,K}$,
and hence the two objective values are equal. We therefore consider
$n<K$ and denote $k=\tilde{\pi}_{n+1}$. The difference between the objective values of
problems~\eqref{eqn:un_primal_sep_fixed_budget_b_abs} and
\eqref{eqn:un_primal_sep} is 
\begin{align*}
|{\bf b}^i(t)|^T \int_{R_{i,n}} (\tilde{\bepsilon}(s)- \hat{\bepsilon}(s)) ds = & |{ b}^i_k(t)| \int_{R_{i,n}}(\tilde{\epsilon}_k(s)- \hat{\epsilon}_{k}(s)) ds\\
&= |{b}_{k}^i(t)| \left(\int_{R_{i,n}} (\Gamma(s) - n) ds - \int_{R_{i,n}} (\Gamma_{i,n} -n) ds \right)\\
&= |{ b}^i_k(t)|\left(\int_{R_{i,n}} \Gamma(s) ds - \int_{R_{i,n}} \Gamma(s) ds \right)= 0,
\end{align*}
where equalities follow from the definitions of  $k$ and $\Gamma_{i,n}$.
\end{proof}

The next lemma shows that the budgeted uncertainty on $\bepsilon(s)$ can be replaced by a budgeted uncertainty on the normalized integrals of $\bepsilon$ over the time intervals. This converts \eqref{eqn:un_primal_sep_fixed_budget_b_abs} from an optimization over a function on a continuous time domain into one over finitely many variables for each $t$, which can then be dualized by the standard techniques of finite-dimensional robust optimization. Specifically, for any $i\in [p]$ and $n\in\{0,1,\ldots, K\}$, defining the average uncertainty vector as
\[
\bxi_{i,n} \equiv \frac{\int_{R_{i,n}}\bepsilon(s) ds}{S_{i,n}}, 
\]
we prove that there exists an equivalent robust problem in terms of $\bxi_{i,n}$.
Thus, we are able to reformulate the problem using classic dualization techniques used in finite-dimensional robust optimization. 

\begin{lemma}
\label{xi_uncertainty_lemma}
Let $\Xi_{i,n}$ be the set of all possible realizations of $\bxi_{i,n}$, that is $$\Xi_{i,n} \equiv \left\{ \bxi\in \Real^{K}: \bxi= \frac{\int_{R_{i,n}}\bepsilon(s) ds}{S_{i,n}}  ,\; \bepsilon(s) \in \U_{i,n} \right\},$$
and let 
\begin{equation} \label{xi_uncertainty_set} \tilde{\U}_{i,n} \equiv \left\{\bxi \in \Real^K: ||\bxi||_1 \leq \Gamma_{i,n}, 0 \leq \bxi \leq 1 \right\}.
\end{equation}
Then $\Xi_{i,n} = \tilde{\U}_{i,n}$.
\end{lemma}
\begin{proof}
We start by showing that $\Xi_{i,n} \subseteq \tilde{\U}_{i,n}$. Let $\hat{\bxi}= \frac{\int_{R_{i,n}} \hat{\bepsilon}(s)ds}{S_{i,n}}   \in \Xi_{i,n}$ for some realization $\hat{\bepsilon}(s) \in \U_{i,n}$. It holds that 
$$\hat{\xi}_k =   \frac{\int_{R_{i,n}} \hat{\epsilon}_k(s) ds}{S_{i,n}}  \leq \frac{1}{S_{i,n}} \int_{R_{i,n}}  \hat{\epsilon}_k(s)   ds\leq \frac{1}{S_{i,n}} \int_{R_{i,n}} 1  ds = 1 ,$$
where the first inequality follows from the triangle inequality and the second inequality follows from $\hat{\epsilon}_k(s)\le 1$, which holds by the definition $\U_{i,n}$. The nonnegativity of entries in $\hat{\bxi}$ follows from the non-negativity of $\hat{\bepsilon}$. Furthermore, 
\begin{align*}
\norm{\hat{\bxi}}_{1} &= \sum_{k=1}^{K} \hat{\xi}_k 
= \sum_{k=1}^{K}  \frac{\int_{R_{i,n}} \hat{\epsilon}_k(s) \, ds}{S_{i,n}}  %
= \frac{1}{S_{i,n}} \int_{R_{i,n}} \sum_{k=1}^{K} |\hat{\epsilon}_k(s)|  \, ds %
\leq \frac{1}{S_{i,n}} \int_{R_{i,n}} \Gamma_{i,n} \, ds %
= \Gamma_{i,n}.
\end{align*}
where the first equality uses $\hat{\bepsilon}\ge\bb{0}$, and the inequality follows from
the definition of $\U_{i,n}$. Hence, $\hat{\bxi} \in \tilde{\U}_{i,n}$.
On the other hand, let $\tilde{\bxi} \in \tilde{\U}_{i,n}$, and define $\hat{\bepsilon}:R_{i,n}\rightarrow \Real^K$ such that for every $s\in R_{i,n}$ and $k\in[K]$ 
$$  \hat{\epsilon}_k(s) = \tilde{\xi}_k\in[0,1].$$  
By definition,
$\tilde{\bxi}=\frac{\int_{R_{i,n}}\hat{\bepsilon}(s)ds}{S_{i,n}}$,
and $$\norm{\hat{\bepsilon}(s)}_1 = \sum_{k=1}^K \hat{\epsilon}_k(s) = \sum_{k=1}^K [\tilde{\bxi}]_k =  \norm{\tilde{\bxi}}_1 \leq   \Gamma_{i,n}.$$ 
Thus, $\hat{\bepsilon} \in \U_{i,n}$ and therefore $\tilde{{\bxi}}\in \Xi_{i,n}$, implying  $\tilde\U_{i,n} \subseteq \Xi_{i,n}$. %
\end{proof}

Using the result in Lemma \ref{xi_uncertainty_lemma} we are now ready to show a finite-dimensional problem with the same objective as problem \eqref{eqn:un_primal_sep}.

\begin{lemma} \label{lem:dual_problem_for_subp} 
Let $l\in [p], i \in [l-1]$, and $t\in(T_{l-1},T_{l}]$. Then, for all $n \in \{0,\dots, K\}$ the associated problem \eqref{eqn:un_primal_sep} has the same objective value as the following problem

\begin{equation}\label{eq:dual_uncr_set1}
\begin{aligned}[t]
\min_{\bv_{i,n}\ge 0,\; \mu_{i,n}\ge 0}\quad
& \mu_{i,n}\Gamma_{i,n} + \bv_{i,n}^\top \be\\
\text{s.t.}\quad
& S_{i,n}|{\bf b}^i(t)| - \bv_{i,n} - \mu_{i,n}\be \le \bb{0}.
\end{aligned}
\end{equation}
\end{lemma}
\begin{proof}
Using Lemma \ref{xi_uncertainty_lemma} we can replace problem \eqref{eqn:un_primal_sep_fixed_budget_b_abs} associated with $t \in (T_{l-1},T_{l}], i \in [l-1]$ with the following classic LP problem:  
 \begin{align}
\label{eqn:primal_xi}
    \max\quad&   S_{i,n} |{\bf b}^i(t)|^\top \bxi_{i,n}   \\ 
 \sto \quad &   ||\bxi_{i,n}||_{1}  \leq \Gamma_{i,n}, \nonumber\\ 
 \quad & 0\leq \bxi_{i,n} \leq 1. \nonumber  
  \end{align} 
Its dual problem is given by \eqref{eq:dual_uncr_set1}.

Since problem \eqref{eqn:primal_xi} is bounded, and can be transformed into a linear optimization problem, strong duality holds for the problem and its dual, and both have the same objective value. 
\end{proof}

Finally, the next theorem aggregates the results in this subsection to obtain the robust counterpart for inequality~\eqref{eq:general_cons_epsilon_ext} with uncertainty set $\U_1$.

\begin{theorem}
Let $k\in[K]$, $l\in[p]$, and $t\in(T_{l-1},T_l]$.  
Then inequality~\eqref{eq:general_cons_epsilon_ext} with uncertainty set
$\U_1$ defined in~\eqref{def:uncertainty_set_1} is satisfied if and only if
there exist auxiliary variables
$\bv_{i,n}$ and $\mu_{i,n}$,
for all $i\in[l-1]$ and $n\in[K]_0$,
such that the following constraints hold:
\begin{align}
\label{eq:cons_deterministic_1}
&\sum_{i=1}^{l-1}\sum_{n=0}^{K}
\bigl(
\mu_{i,n}\Gamma_{i,n}
+\bv_{i,n}^\top \be
\bigr)
+ \tilde{\lambda}^0_{k}
\int_{T_{l-1}}^{t} \min\{\Gamma(s),1\}\, ds
\;\le\; z(t),\\
\label{eq:cons_dual_feasibility_1}
&S_{i,n}|{\bf b}^i(t)|
-\bv_{i,n}
- \mu_{i,n}\be
\;\le\; \bb{0},
\qquad
\forall\, i\in[l-1],\, n\in[K]_0,\\
&\bv_{i,n}\ge \bb{0},
\;
\mu_{i,n}\ge 0,\qquad\qquad\qquad\quad\;\forall\, i\in[l-1],\, n\in[K]_0.\nonumber
\end{align}
\end{theorem}

\begin{proof}
As discussed above, \eqref{eq:general_cons_epsilon_ext} can be computed as the sum of the optimal values of problem \eqref{eqn:seperable_prob_i} with $i \in [l-1]$ and \eqref{eqn:seperable_prob_t}. The last term in left-hand side of \eqref{eq:cons_deterministic_1} is due to the optimal value of \eqref{eqn:seperable_prob_t} discussed above. %
Moreover, for each $i\in[l-1]$, problem \eqref{eq:seperable_prob_i_abs} decomposes into the separable problems~\eqref{eqn:un_primal_sep} over $n\in\{0,\ldots,K\}$, and by Lemma~\ref{lem:dual_problem_for_subp} each of these has the same optimal value as \eqref{eq:dual_uncr_set1}. Substituting accordingly yields the double sum in \eqref{eq:cons_deterministic_1}. Finally, as is standard in reformulations of robust optimization problems, the minimum over the auxiliary dual variables $\mu_{i,n}$ and $\bv_{i,n}$ can be removed and replaced by feasibility constraints. 
\end{proof}

\subsection{Budgeted Uncertainty over Time Intervals} \label{sec:uncert_over_time}

The next budgeted uncertainty set we consider limits the accumulated relative perturbations of arrival rates for each buffer over predefined time intervals. Mathematically, we describe this uncertainty set as:
\begin{equation}\label{def:uncertainty_set_2}
\U_2
    = \left\{
        \bepsilon(\cdot)\in \M^{[0,T],K} \;\middle|\;
        \|\bepsilon(s)\|_{\infty} \le 1 \; \forall s\in[0,T],\;
        \int_{T_{l-1}}^{T_l} |\bepsilon(s)|\,ds \le \tilde{\bGamma}_l,\;
        \forall\, l\in[p]
      \right\},
\end{equation}
where  $\tilde{\bGamma}_l \in {\Real}^K_+$, and without loss of generality we assume that it satisfies $\tilde{\bGamma}_l \leq T_{l}-T_{l-1} ,\forall l \in [p]$. 

As in the previous subsection, our goal is to obtain a robust 
counterpart of constraint~\eqref{eq:general_cons_epsilon_ext} for a fixed $t$, 
with a finite number of constraints and variables. This stands in contrast to the 
current formulation, which involves infinitely many variables and constraints 
induced by the uncertainty set $\U_2$.

We start by defining the accumulated uncertainty vector $\bxi_{t}^{i}  \in \Real^{K}$ for each $i\in[p]$ and $t\in (T_{i-1},T_i]$ as  
\begin{equation}\label{eq:def_xi_vertical}
\bxi_{t}^{i} = \int_{T_{i-1}}^{t} \bepsilon(s) ds.
\end{equation}
In the following lemma, we prove that the integral operation over the time intervals $\{(T_{i-1},T_i]\}_{i=1}^{l-1} \cup (T_{l-1},t]$ for fixed $t \in (T_{l-1},T_l]$  applied on the uncertain parameter $\bepsilon$ transforms the uncertainty set $\U_2$ into a box-type uncertainty set over the concatenated vector $\bxi=(\bxi^1_{T_1},\bxi^2_{T_2},\ldots,\bxi^{l-1}_{T_{l-1}},\bxi^{l}_t)$.

\begin{lemma} 
\label{lem:xi_uncertainty_lemma_second}
For some $t \in (T_{l-1},T_l], \quad l \in [p]$, let $\Xi_{t}$ be the set of all possible realizations of $\{\bxi_{T_{i}}^{i}\}_{i\in[l-1]}$ and $\bxi_{t}^{l}$%
attainable under $\U_2$, that is
$$\Xi_{t} \equiv\left\{ (\bxi^{1},\bxi^{2},\ldots,\bxi^{l})\in\Real^{Kl}\;\left|\;
\begin{array}{ll} \exists \bepsilon \in \U_2,\;
    \bxi^{i} =  \int_{T_{i-1}}^{T_i} \bepsilon(s) ds, &\forall i \in [l-1],  \\
    \bxi^{l} = \int_{T_{l-1}}^{t} \bepsilon(s) ds, & 
\end{array} \right.\right\},$$ and let

\begin{equation}
\tilde{\U}_{t} \equiv  \left\{ (\bxi^1,\bxi^2,\ldots,\bxi^l) \in  \Real^{Kl} \;|\; |\bb{\xi^i}|  \leq  \tilde{\bGamma}_i,\;   \forall i \in [l-1],\; |\bxi^l| \leq \min \{t-T_{l-1},\tilde{\bGamma}_l\}  \right \}. 
\end{equation}
Then $\Xi_{t} = \tilde{\U}_{t}$.
\end{lemma}
\begin{proof}
We start by showing that $\Xi_{t} \subseteq \tilde{\U}_{t}$. Let $\hat{\bxi} \in \Xi_t$. It holds that there exist some realization $\hat{\bepsilon}(s) \in \U_2$, for which 
\begin{align*}|\hat{\bxi}_{t}^l| &= \left|\int_{T_{l-1}}^{t} \hat{\bepsilon}(s) ds\right| \leq \int_{T_{l-1}}^{t} |\hat{\bepsilon}(s)| ds %
\leq \min\{t-T_{l-1}, \tilde{\bGamma}_l\} , \\
|\hat{\bxi}_{T_i}^{i}| &= \left|\int_{T_{i-1}}^{T_i} \hat{\bepsilon}(s) ds\right| \leq \int_{T_{i-1}}^{T_i} |\hat{\bepsilon}(s)| ds%
\leq \tilde{\bGamma}_i, \quad \forall i \in [l-1],
\end{align*}
where the first inequality in both equations follows from the triangle inequality and the second inequality follows from the definition of $\U_2$, thus, $\hat{\bxi} \in \tilde{\U}_{t}$.

On the other hand, let $\tilde{\bxi}\in \tilde{\U}_t$.
Define the realization $\hat{\bepsilon}(s)$ 
 as:
\[ \hat{\bepsilon}(s) = \begin{cases} \frac{\tilde{\bxi}^{i}}{T_i-T_{i-1}}, & s\in(T_{i-1},T_i], \quad i \in [l-1],\\
\frac{\tilde{\bxi}^{l}}{t-T_{l-1}}, &  s\in (T_{l-1} ,t],   \\ 0 & else.
\end{cases} \] 
Since $\tilde{\bGamma}_i \leq T_{i}-T_{i-1}$ for $i\in[l-1]$ and $|\tilde{\bxi}^l|\le t-T_{l-1}$ by the definition of $\tilde{\U}_t$, we conclude that $\norm{\hat{\bepsilon}(s)}_{\infty} \leq 1$ for all $s\in [0,T]$.
In addition, $$\int_{T_{l-1}}^{T_l} |\hat{\bepsilon}(s)| ds = \int_{T_{l-1}}^{t} |\hat{\bepsilon}(s)| ds 
 = \int_{T_{l-1}}^{t} \frac{|\tilde{\bxi}^{l}|}{t-T_{l-1}} ds = |\tilde{\bxi}^{l}| \leq \tilde{\bGamma}_l,$$
$$\int_{T_{i-1}}^{T_i} |\hat{\bepsilon}_k(s)| ds 
 = \int_{T_{i-1}}^{T_i} \frac{|\tilde{\bxi}^{i}|}{T_i-T_{i-1}} ds = |\tilde{\bxi}^{i}| \leq \tilde{\bGamma}_i,\quad \forall i \in [l-1],$$ 
 where the last inequalities follow from the definition of $\tilde{\U}_t$. Thus, we found a realization $\hat{\bepsilon}(s) \in \U_2$ such that  $\tilde{\bxi} = (\tilde{\bxi}_{T_{1}}^{1},\tilde{\bxi}_{T_{2}}^{2},\dots, \tilde{\bxi}_{T_{l-1}}^{l-1},\tilde{\bxi}_{t}^{l})$, is an element of $\Xi_t$ associated with $\hat{\bepsilon}(s)$,  establishing that the two sets are the same.
  Together with the first inclusion, this establishes $\Xi_t=\tilde{\U}_t$.
\end{proof}

Using this result,
we turn to reformulate constraint \eqref{eq:general_cons_epsilon_ext}  for uncertainty
set $\U_2$.

\begin{theorem}
    Let $k\in[K]$, $l\in [p]$, and $t\in(T_{l-1},T_{l}]$. Then inequality~\eqref{eq:general_cons_epsilon_ext} with uncertainty set $\U_2$ given by \eqref{def:uncertainty_set_2} is equivalent to 

\begin{equation}\label{eq:RO_for_second_uncertainty}
 \sum_{i=1}^{l-1}  |{\bf b}^i(t)|^\top \tilde{\bGamma}_i + \tilde{\lambda}^0_{k} \Omega_k(t) \leq z(t), 
\end{equation}

where $\Omega_k(t)= \min \left\{t-T_{l-1}, [\tilde{\bGamma}_l]_{k} \right\}$.
\end{theorem}
\begin{proof}
Since \eqref{eq:general_cons_epsilon_ext} involves only terms of the form $(\bxi^1_{T_1},\bxi^2_{T_2},\ldots,\bxi^{l-1}_{T_{l-1}},\bxi^{l}_t)$, as defined above, one can replace the maximization over $\U_2$ with the maximization over $\Xi_t$. Moreover, 
applying Lemma \ref{lem:xi_uncertainty_lemma_second} to constraint \eqref{eq:general_cons_epsilon_ext} associated with $t \in (T_{l-1},T_{l}], l \in [p]$ we can replace $\Xi_t$ by $\tilde{U}_t$ yielding the following equivalent constraint 
\begin{equation} \label{eq:cons_2_after_lemma}
\max_{\bxi\in \tilde{\U}_t}\sum_{i=1}^{l-1} {\bf b}^i(t)^\top \bxi^i -    \tilde{\lambda}^{0}_k {\xi}^l_k\leq z(t).
\end{equation}
Moreover, using the separability of $\tilde{\U}_t$ over $i\in[l]$ and $k\in[K]$, we have
$$\max_{\bxi\in \tilde{\U}_t}\sum_{i=1}^{l-1} {\bf b}^i(t)^\top \bxi^i -    \tilde{\lambda}^{0}_k {\xi}^l_k=\sum_{i=1}^{l-1} \max_{|\bxi^i|\leq \tilde{\bGamma}_i}{\bf b}^i(t)^\top \bxi^i+\max_{|{\xi}^l_k|\leq \min\{t-T_{l-1},[\tilde{\bGamma}_l]_k\}}-\tilde{\lambda}^{0}_k {\xi}^l_k.$$
which equals the left-hand side of \eqref{eq:RO_for_second_uncertainty}. %
\end{proof}

\subsection{Budgeted Uncertainty over Buffers and Time Intervals}
\label{sec:uncertainty_set_3_combined} 

The uncertainty set considered in this subsection is the amalgamation of the two uncertainty sets $\U_1$ and $\U_2$. On the one hand, we limit the uncertainty in arrival rate to each buffer over predefined time intervals using the budget vectors $\{\tilde{\bGamma}_l\}_{l\in[p]} \subseteq \mathbf{R}^K_+$, and on the other hand, we bound the uncertainty in the sum of arrival rates over all buffers at each specific point in time using the budget function $\Gamma(\cdot)$. %
Similarly to sections \ref{sec:uncert_in_all_buffers} and \ref{sec:uncert_over_time}, we assume throughout this subsection that $\Gamma(\cdot)$ satisfies Assumption~\ref{ass:Gamma} and $\tilde{\bGamma}_l\leq T_{l}-T_{l-1}$.

Mathematically, the uncertainty set is given by: 
\begin{equation} \label{def:uncertainty_set_3}
\begin{aligned}
\U_3 
= \left\{
    \bepsilon(\cdot)\in\M^{[0,T],K} 
    \;\middle|\;
   \begin{array}{l}  \quad\|\bepsilon(t)\|_\infty \le 1,\;  \qquad\;\;\forall t \in [0,T],  \\[4pt]
    \quad \int_{T_{i-1}}^{T_i} |\bepsilon(s)| \, ds \le \tilde{\bGamma}_i,\; \forall i \in [p],  \\[4pt]
    \quad \|\bepsilon(t)\|_1 \le \Gamma(t),\;  \qquad\forall t \in [0,T]
    \end{array}
\right\}.
\end{aligned}
\end{equation}

Note that this uncertainty set is a generalization of $\U_1$ and $\U_2$ as setting $\tilde{\bGamma}_i = (T_i-T_{i-1})\be$ for each $i \in [p]$ results in uncertainty set $\mathcal{U}_1$, and setting $\Gamma(t) = K$ for all $t\in[0,T]$ results in uncertainty set $\U_2$.

The reformulation of constraint \eqref{eq:general_cons_epsilon_ext} for uncertainty set $\U_3$ is given in the following theorem, where we reuse the notations 
 $R_{i,n}$, $S_{i,n}$, $\Gamma_{i,n}$, and  $\bxi_{i,n}$ defined in Section~\ref{sec:uncert_in_all_buffers}.

\begin{theorem}
\label{thm:robust_counterpart_constraint_in_3}
Let $k\in[K]$, $l\in[p]$, and $t\in(T_{l-1},T_l]$.  
Then inequality~\eqref{eq:general_cons_epsilon_ext} with uncertainty set
$\U_3$ defined in~\eqref{def:uncertainty_set_3} is satisfied if and only if
there exist additional variables $\bgamma_i\in \Real^K_+$, $\bv_{i,n}\in\Real^K_+$ and $\mu_{i,n}\in \Real_+$ for all $i\in[l-1]$, $n\in[K]_0$
such that the following  inequalities hold:
\begin{subequations}
\label{eq:cons_deterministic_3_linear}
\begin{align}
\label{eq:cons_deterministic_3_linear_a}
\sum_{i=1}^{l-1} \left[\sum_{n=0}^{K}
\Bigl(
\mu_{i,n}\Gamma_{i,n}
+\bv_{i,n}^\top \be\Bigr) 
+ \bgamma_i^\top \tilde{\bGamma}_i \right]
+ \tilde{\lambda}^0_{k}
 \Omega_k(t)
&\le z(t), \\[0.5em]
\label{eq:cons_deterministic_3_linear_b}
S_{i,n}\bigl(|{\bf b}^i(t)|-\bgamma_i\bigr)
-\bv_{i,n} - \mu_{i,n}\be
&\le \bb{0},
\;
\forall\, i\in[l-1],\ \forall\, n\in[K]_0,
\end{align}
\end{subequations}
where $\Omega_k(t)=\min\!\left\{
\int_{T_{l-1}}^{t}\min\{\Gamma(s),1\}\, ds,\,
[\tilde{\bGamma}_l]_{k}
\right\}$.
\end{theorem}
\begin{proof}To find the deterministic robust counterpart for constraint \eqref{eq:general_cons_epsilon_ext} associated with $l \in [p], t \in (T_{l-1},T_l]$ we need to maximize the left-hand side of the constraint under uncertainty set $\U_3$ defined in \eqref{def:uncertainty_set_3}.
As in subsection \ref{sec:uncert_in_all_buffers},
due to separability of the constraint
and uncertainty set over the time intervals $\{(T_{l-1},T_l]\}_{l=1}^{p}$,
we decompose the maximization problem over $\U_3$ %
into $l$ subproblems. The $l$th subproblem is given by
\begin{align} \label{eq:second_type_subpro_3}
     \max \quad &  - \int_{T_{l-1}}^t \tilde{\lambda}^0_{k} \epsilon_{k}(s) ds  \\
s.t. \quad &   \sum_{k=1}^K |\epsilon_k(s)|  \leq \Gamma(s), \quad \forall s \in (T_{l-1},T_l],\nonumber \\ 
 \quad &  |\bepsilon(s)|  \leq 1, \quad  \forall s \in (T_{l-1},T_l], \nonumber \\ \quad & \int_{T_{l-1}}^{T_l} |\bepsilon(s)| ds \leq \tilde{\bGamma}_l,\nonumber
\end{align}
which has the trivial optimal value $\tilde{\lambda}_k^{0}\Omega_k(t)$, %
obtained by setting $\epsilon_{k'}(s)=0$ for all $k'\in[K]\setminus\{k\}$ and taking $\epsilon_k(s)$ at its minimum possible value under the pointwise constraints for $s\in(T_{l-1},t]$ until the interval budget is reached.
The remaining $l-1$ subproblems, indexed by $i\in[l-1]$, are given by
\begin{align}
\label{eqn:un_primal_3}
     \max  \quad &{{\bf b}}^i(t)^\top \int_{T_{i-1}}^{T_{i}} \bepsilon(s) ds  \\ 
s.t. \quad &   \sum_{k=1}^K |\epsilon_k(s)|  \leq \Gamma(s), \quad \forall s \in (T_{i-1},T_i], \nonumber \\ 
 \quad & |\bepsilon(s)|  \leq 1, \quad  \forall s \in (T_{i-1},T_i], \nonumber 
 \\ \quad & \int_{T_{i-1}}^{T_{i}} |\bepsilon(s)| ds \leq \tilde{\bGamma}_i. \nonumber
\end{align}

From symmetry of the constraints around zero,
we note that the sign of $\bepsilon(s)$ will always correspond to the sign of  $\bbb^i(t)$.
Therefore, problem \eqref{eqn:un_primal_3} is equivalent to the following problem
\begin{subequations}\label{eqn:un_primal_3_without_abs_val}
    \begin{align}
        \max_{\epsilon(\cdot)}  \quad & |{\bbb}^i(t)|^\top \int_{T_{i-1}}^{T_{i}} \bepsilon(s) ds \nonumber \\ 
        \text{s.t.} \quad 
        & \sum_{k=1}^K \epsilon_k(s)  \leq \Gamma(s), \quad \forall s \in (T_{i-1},T_i], \label{constraint1a} \\ 
        & 0 \leq \bepsilon(s)  \leq 1, \quad \forall s \in (T_{i-1},T_i], \label{constraint1b} \\
        & \int_{T_{i-1}}^{T_{i}} \bepsilon(s) ds \leq \tilde{\bGamma}_i. \label{constraint1c}
    \end{align}
\end{subequations}

Any feasible solution $\bepsilon(\cdot)$ to problem \eqref{eqn:un_primal_3_without_abs_val} can be adapted to a feasible solution $\operatorname{sign}({\bf b}^i(t)) \bepsilon(\cdot)$ to problem \eqref{eqn:un_primal_3}, and vice versa, $|\bepsilon(\cdot)|$ is feasible to problem \eqref{eqn:un_primal_3_without_abs_val}, for any feasible solution $\bepsilon(\cdot)$ to \eqref{eqn:un_primal_3}. Since $b^i_k(t)\epsilon_k(s)\leq |{ b}^i_k(t)| |\epsilon_k(s)|$, there must be an optimal solution of \eqref{eqn:un_primal_3} satisfying $\sign(\epsilon_k(s))=\sign(b^i_k(t))$, and both problems have the same optimal value. %

Note that \eqref{eqn:un_primal_3_without_abs_val} is an SCLP problem, and dualizing it only with respect to the last constraint results in
\begin{equation*} %
\begin{aligned}[t]
\min_{\bgamma_i \geq 0} \max_{\bepsilon(\cdot)} & \quad \left(|{\bf b}^i(t)| - \bgamma_i\right)^\top \int_{T_{i-1}}^{T_{i}} \bepsilon(s) \, ds + \bgamma_i^\top \tilde{\bGamma}_i \\
\text{s.t.} & \quad \sum_{k=1}^K \epsilon_k(s) \leq \Gamma(s),  \quad && \forall s \in (T_{i-1},T_i], \\
& \quad 0 \leq \bepsilon(s) \leq 1,  \quad && \forall s \in (T_{i-1},T_i].
\end{aligned}
\end{equation*}

{Viewing problem~\eqref{eqn:un_primal_3_without_abs_val} over the space of bounded
measurable functions equipped with the supremum norm, Slater's condition is satisfied, and strong duality holds \cite[Theorem 2.9.2]{zalinescu2002convex}}.
Looking at the inner maximization problem, the separability of the constraints and objective over time allows to decompose the problem into $K+1$ subproblems associated with the restricted domains $R_{i,n}$ for $n\in[K]_0$, as in Section \ref{sec:uncert_in_all_buffers}. Applying  the dualization argument of Lemma \ref{lem:dual_problem_for_subp} to each subproblem associated with $R_{i,n}$, for $n\in[K]_0$, yields
the following problem: 
\begin{equation}
\begin{aligned}[t]\label{eq:final_dual_RO_for_subpr_3}
\min_{\bgamma_i\ge 0,\;\bv_{i,n}\ge 0,\; \mu_{i,n}\ge 0} \quad \sum_{n=0}^{K} 
& \left( \mu_{i,n}\Gamma_{i,n} +\bv_{i,n}^\top \be \right) + \bgamma_i^\top \tilde{\bGamma}_i  \\
\text{s.t.}\quad
& S_{i,n}(|{\bf b}^i(t)|- \bgamma_i) -\bv_{i,n} - \mu_{i,n}\be \le \bb{0}, \, n \in[K]_0.
\end{aligned}
\end{equation}

Substituting the objective terms of the dual problems
\eqref{eq:final_dual_RO_for_subpr_3} for all relevant $i\in[l-1]$
back into~\eqref{eq:general_cons_epsilon_ext}, together with the robust solution
of the subproblem~\eqref{eq:second_type_subpro_3},
yields the robust counterpart of
constraint~\eqref{eq:general_cons_epsilon_ext}, given by \eqref{eq:cons_deterministic_3_linear} with additional nonnegativity constraints on the variables.

\end{proof}

\subsection{Reformulation of the Robust Effort Model}

As mentioned above, since uncertainty set $\U_3$ is a generalization of both $\U_1$ and $\U_2$, we only provide the  reformulation of the constraints of problem \eqref{prob:server_effort_model_only_uncertainty_general} for this set. Applying Theorem \ref{thm:robust_counterpart_constraint_in_3} to all uncertain constraints of problem \eqref{prob:server_effort_model_only_uncertainty_general}, we note that since the constraints need to be satisfied for all $t\in[0,T]$, the dual variables become a function of $t$.

Thus, defining  $\chi_\ell(t):=\mathbf{1}_{\{t> T_{\ell}\}}$ for each $\ell\in[p-1]$, and using the transformation $\Phi^{\ell}(t)=\Phi^{\ell+}(t)-\Phi^{\ell-}(t)$, where  $\Phi^{\ell+}(t),\Phi^{\ell-}(t)\geq 0$ for all $t\in[0,T]$,
the following formulation is equivalent to problem
\eqref{prob:server_effort_model_only_uncertainty_general} with $\U=\U_3$.
\begin{subequations}\label{eq:final_linear_RC}
\begin{align}
\min
\ \max_{\bepsilon \in \U}
\quad &
F\!\left(\tilde{\bleta}^{0}, \{\tilde{\bPhi}^{q}\}_{q=1}^{p-1}, \bepsilon\right)
\label{eq:final_linear_RC:a}
\\[0.8em]
\text{s.t.}\quad 
& 
\bG_{k}
\int_0^t \tilde{\bleta}^0(\tau)\, d\tau
+
\sum_{\ell=1}^{p-1}
\Biggl[
\sum_{n=0}^{K}
\Bigl(
\Gamma_{\ell,n}\mu^{1}_{k,\ell,n}(t)
+
\be^\top\bv^{1}_{k,\ell,n}(t)
\Bigr)
+
\tilde{\bGamma}_\ell^\top\bgamma^{1}_{k,\ell}(t)
\Biggr]\nonumber\\
&\quad
\le
\alpha_{k} + \hat{\lambda}_{k} t - \Omega_{k}(t)\tilde{\lambda}_{k},
\hspace{8em}\quad\;\;
k \in [K],\,
t \in [0,T],
\label{eq:final_linear_RC:b}
\\[1em]
& 
S_{\ell,n}
\left(\blbeta^1_{k,\ell}(t)
-
{\bgamma}^{1}_{k,\ell}(t)
\right)
-
\bv^{1}_{k,\ell,n}(t)
-
\be\,\mu^{1}_{k,\ell,n}(t)
\le 0,
\nonumber\\%
&\hspace{11em}\quad
\;\, k \in [K],\,
\ell \in [p-1],\,
n \in [K]_0,\,
t \in [0,T],\nonumber
\\[1em]
&\left|
\int_{0}^{t}
\left[\tilde{\bPhi}^{\ell+}(\tau)-\tilde{\bPhi}^{\ell-}(\tau)\right]^\top \bG_{k}^\top\, d\tau
-
\chi_\ell(t)\be_k\tilde{\lambda}_k
\right|-\blbeta^1_{k,\ell}(t)\leq 0\nonumber\\
&\hspace{15em}
\forall\, k \in [K],\,
\ell \in [p-1],
\,t \in [0,T],\nonumber%
\\[1em]
& \sum_{j:\,\tilde{s}(j)=i}\tilde{\eta}^{0}_j(t)+\sum_{\ell=1}^{p-1} \Biggl[ \sum_{n=0}^{K} \Bigl( \Gamma_{\ell,n}\mu^{2}_{i,\ell,n}(t) + \be^\top\bv^{2}_{i,\ell,n}(t) \Bigr) + \tilde{\bGamma}_{\ell}^\top \bgamma^{2}_{i,\ell}(t) \Biggr] \le 1 , \nonumber\\
&\;\hspace{20em}
\,\; i \in [I],\,
t \in [0,T],
\label{eq:final_linear_RC:c}
\\[1em]
& S_{\ell,n} \Bigl(  \blbeta^2_{i,\ell}(t) - \bgamma^{2}_{i,\ell}(t) \Bigr) - \bv^{2}_{i,\ell,n}(t) - \be\,\mu^{2}_{i,\ell,n}(t) \le 0, \nonumber\\ &\hspace{11em} \,\, i \in [I],\, \ell \in [p-1],\, n \in [K]_0,\, t \in [0,T],\nonumber%
\\[1em]
&\left|\sum_{j:\,\tilde{s}(j)=i} [\tilde{\bPhi}_{j}^{\ell+}(t)-\tilde{\bPhi}_{j}^{\ell-}(t)]^\top\right|-{\blbeta}^2_{i,\ell}(t)\leq 0,\nonumber\\[-1em]
 &\hspace{15em} \,\,
\,\, i \in [I],\,
\ell \in [p-1],
\,t \in [0,T],\nonumber\\[1em]
& \tilde{\eta}_{j}^{0}(t) - \sum_{\ell=1}^{p-1} \Biggl[ \sum_{n=0}^{K} \Bigl( \Gamma_{\ell,n}\mu^{3}_{j,\ell,n}(t) + \be^\top\bv^{3}_{j,\ell,n}(t) \Bigr) + \tilde{\bGamma}_\ell^\top\bgamma^{3}_{j,\ell}(t) \Biggr] \ge 0, \nonumber\\ &\hspace{20em}\, \,\,j \in [J],\, t \in [0,T],
\label{eq:final_linear_RC:d}
\\[1em]
&S_{\ell,n} \Bigl( 
[\tilde{\bPhi}_{j}^{\ell+}(t)+\tilde{\bPhi}_{j}^{\ell-}(t)]^\top- \bgamma^{3}_{j,\ell}(t) \Bigr) - \bv^{3}_{j,\ell,n}(t) - \be\,\mu^{3}_{j,\ell,n}(t) \le 0, \nonumber\\ &\hspace{11em}\, \, j \in [J],\, \ell \in [p-1],\, n \in [K]_0,\, t \in [0,T],\nonumber
\\[1em]
&
\tilde{\bPhi}^{\ell+}(t) = \tilde{\bPhi}^{\ell-}(t) = 0,\, \blbeta^r_{\cdot,\ell}(t)=0,\;
 r\in[2],\,\ell \in [p-1],\,t \in [0,T_\ell],
\label{eq:final_linear_RC:e}
\\[1em]
&
\mu^{r}_{\cdot,\ell,\cdot}(t)=0,
\,
\bv^{r}_{\cdot,\ell,\cdot}(t)=
\bgamma^{r}_{\cdot,\ell}(t)=0,
\,\,\, r\in[3],\,
\ell\in[p-1],\,
t\in[0,T_\ell],\label{eq:final_linear_RC:f_gamma}
\\[1em]
&\blbeta^r_{\cdot,\ell}(t) \ge 0,\, \tilde{\bPhi}^{\ell+}(t) \geq 0,\,
\tilde{\bPhi}^{\ell-}(t)\geq 0,\;r \in[2],\ell\in[p-1],\, t\in[0,T].\nonumber\\[1em]%
&
\bv^{r}_{\cdot,\cdot,\cdot}(t) \ge 0,\,
\mu^{r}_{\cdot,\cdot,\cdot}(t) \ge 0,\,
\bgamma^{r}_{\cdot,\cdot}(t) \ge 0, \hspace{5em}r \in[3],\,t\in[0,T]
,
\nonumber
\end{align}
\end{subequations}

The role of the activation constraints
\eqref{eq:final_linear_RC:e}--\eqref{eq:final_linear_RC:f_gamma}
is to remove the dependence on the interval index \(l\). Instead of writing a
different constraint on every interval \((T_{\ell-1},T_\ell]\), all constraints are now
imposed on the full horizon \(t\in[0,T]\). Variables associated with segment
\(\ell\) are forced to be zero before their activation time \(T_\ell\), and therefore
they do not affect the constraints before that time. Moreover, constraints \eqref{eq:final_linear_RC:f_gamma} are in fact redundant due to constraint \eqref{eq:final_linear_RC:e} %
and the use of activation function $\chi_l$.

Note that this problem is equivalent to the SCLP problem \eqref{eq:extended_sclp_form}.   Specifically, the controller $\bu$ is given by the variables $\tilde{\eta}^0,\Phi^{\ell}, {\mu}^{r}_{\cdot,\ell,\cdot},\bv^r_{\cdot,\ell,\cdot},\bgamma^{r}_{\cdot,\ell},\blbeta^2_{\cdot,\ell}$ for all $\ell\in[p-1]$ and $r\in\{2,3\}$, and the state $\bx$ is given by ${\mu}^{1}_{\cdot,\cdot,\cdot},\bv^1_{\cdot,\cdot,\cdot},\bgamma^{1}_{\cdot,\cdot},\blbeta^{1}_{\cdot,\cdot}$.
Moreover, by
Assumption~\ref{ass:Gamma}, the uncertainty-budget function
is piecewise constant, so the integrated budget term
\(\Omega_k(t)\) given in Theorem~\ref{thm:robust_counterpart_constraint_in_3} is piecewise affine. Separately, constraint \eqref{eq:final_linear_RC:d} together with the nonnegativity of the dual variables implies $\tilde{\eta}^0\ge 0$.

\section{Establishing an SCLP Min-Max Algorithm}\label{sec:Algorithm}

Although the constraints in \eqref{eq:final_linear_RC} can be handled directly under the assumed extension of the SCLP--Simplex algorithm~\cite{Weiss2008,Shindin2021} discussed in Section~\ref{sec:preliminaries}, the objective function still poses a significant challenge.
Applying the same reasoning as in Section~\ref{sec:constraints}, the robust counterpart of the objective 
contains the absolute value of an integral over the interval $[0,T]$ involving the problem variables. 
The standard technique of introducing auxiliary variables to replace the absolute value does not 
resolve this issue, since the resulting constraints would themselves involve integrals over the complete interval
$[0,T]$, which are not supported by the SCLP--Simplex algorithm.
Thus, we seek another approach to solve this problem. 

Before discussing algorithms for solving this problem, we establish some desirable properties of the feasible controller set of~\eqref{eq:final_linear_RC} {(i.e., the projection of the feasible solution set onto the controller parameters $\tilde{\eta}$ and $\{\tilde{\Phi}^i\}_{i\in[p-1]}$)}, denoted by $\W$, the uncertainty set $\U$, and the objective function $F$.
We begin by providing two equivalent representations of $F$ that make its dependence on each of its arguments explicit. 

We use the shorthand $\bw=(\tilde{\bleta}^{0},\{\tilde{\bPhi}^{i}\}_{i\in[p-1]})$ for the controller parameters.
Decomposing the rightmost integral in \eqref{eqn:obj_term_raw}, $F$ can be represented as
\begin{equation}\label{eqn:objective_first_form}
\begin{aligned}[t]
F(\bw,\bepsilon)
={}&
\tilde{F}(\tilde{\bleta}^{0},\{\tilde{\bPhi}^{i}\},\bepsilon) +
[\tilde{\blambda}\circ\bc]^\top
\sum_{i=1}^{p}
\int_{T_{i-1}}^{T_i}
\int_{T_{i-1}}^t
\bepsilon(s)\,ds\,dt,
\\
\tilde{F}(\bw,\bepsilon)
={}&
\int_0^T
\bh(t)^\top\tilde{\bleta}^{0}(t)\,dt
\\
&\hspace{-20pt}+
\sum_{i=1}^{p-1}
\left[
\int_{T_{i-1}}^{T_i}\bepsilon(s)\,ds
\right]^\top
\left[
\int_0^T
\tilde{\bPhi}^{i}(t)^\top\bh(t)\,dt
+
(T-T_i)(\tilde{\blambda}\circ\bc)
\right].
\end{aligned}
\end{equation}
Here, $\bh(t)$ is the affine function given by
\begin{equation}\label{eq:h_def}
\bh(t)
:=
-(T-t)\bG^\top\bc
+
(\bd\circ\bmu),
\qquad t\in[0,T].
\end{equation}
Thus, for every fixed $\bepsilon\in\U$, the coefficient functions of
$F(\cdot,\bepsilon)$ are affine in time. 

Alternatively, using Fubini's theorem, the objective can be written as
\begin{equation}\label{eq:objective_first_form}
\begin{aligned}[b]
F(\bw, \boldsymbol{\epsilon})
=\;&
\int_0^T
\bh(t)^\top
\tilde{\boldsymbol{\eta}}^{0}(t)\,dt \\
&\; + \sum_{i=1}^{p-1}
\left(
\int_0^T \bh(t)^\top \tilde{\mathbf{\Phi}}^{i}(t)\,dt
\right)
\left(
\int_{T_{i-1}}^{T_i}\boldsymbol{\epsilon}(t)\,dt
\right) \\
&\; + \int_0^T (T-t)(\tilde{\boldsymbol{\lambda}} \circ \mathbf{c})^\top
\boldsymbol{\epsilon}(t)\,dt \\
=\;&
\int_0^T
\bh(t)^\top
\tilde{\boldsymbol{\eta}}^{0}(t)\,dt
+ \int_0^T
\left[
\bbf(t) + (T-t)(\tilde{\boldsymbol{\lambda}} \circ \mathbf{c})
\right]^\top
\boldsymbol{\epsilon}(t)\,dt ,
\end{aligned}
\end{equation}
where $\bh$ is defined in \eqref{eq:h_def}
and \(\bbf(t)\) is the piecewise constant function defined by
\[
\bbf(t):=\begin{cases}
\int\limits_0^T (\tilde{\mathbf{\Phi}}^{i}(s))^\top  \bh(s)\,ds,
&
t\in (T_{i-1},T_i],\quad i=1,\dots,p-1,\\
\bb{0}, &\text{otherwise}.
\end{cases}
\]
Thus, for fixed $\bw\in\W$, the coefficient function of
$F(\bw,\cdot)$ is piecewise affine in time.

Having characterized the dependence of $F$ on each of its arguments, we next establish boundedness properties of the feasible controller set $\W$.
Without loss of generality,  
we impose  $\tilde{\bPhi}^{i}_{\cdot,l}(t)=0,$ for all $t\in[0,T]$, or equivalently omit this column from the model, whenever the set $\U$ implies that $\bepsilon_l(t)=0$ for almost every $t\in(T_{i-1},T_{i}]$.

\begin{lemma}
    Let $\bw=(\tilde{\bleta}^{0}(\cdot),\{\tilde{\bPhi}^{i}(\cdot)\}_{i\in[p-1]})$ be a feasible controller to problem \eqref{eq:final_linear_RC}. Then, $\tilde{\eta}_{j}^0$ and $\tilde{\Phi}^{i}_{jk}$ for all $j\in [J], k\in [K]$, and $i\in[p-1]$, are uniformly bounded in $L^\infty[0,T]$ and $L^2[0,T]$ over the feasible controller set $\W$. 
\end{lemma}
\begin{proof}
Since the constraints of \eqref{eq:final_linear_RC} are equivalent to those of \eqref{prob:server_effort_model_only_uncertainty_general}, it suffices to establish the result using the latter formulation.
Specifically, since \eqref{prob:server_effort_model_only_uncertainty_general:d} has to apply for all $\epsilon\in \U_3$ and $\epsilon=0$ in particular, 
$\tilde{\bleta}^0(t)\geq 0$. 
Combined with \eqref{prob:server_effort_model_only_uncertainty_general:c} we have that $0\leq \tilde{\eta}_j^0(t)\leq 1$ for all $t\in[0,T]$ and $j\in[J]$. Thus, $\tilde{\eta}_j^0\in L^{\infty}[0,T]$ for all $j\in[J]$, and since $T$ is finite, they are also bounded in $L^{2}[0,T]$.  
To show the same for $\tilde{\bPhi}^i$ we focus on $t\in (T_i,T]$, since it is defined as zero for $[0,T_i]$. Let $l\in [K]$, 
setting ${\bepsilon}(s)=\be_l\alpha_{i,l}(s)$ where $\alpha_{i,l}(s)=\min\{1,\frac{[\tilde{\bGamma}_i]_l}{T_i-T_{i-1}},\Gamma(s)\}$ for $s\in (T_{i-1},T_i]$ and zero otherwise, $\pm\bepsilon\in \U_3$. Plugging both into \eqref{prob:server_effort_model_only_uncertainty_general:d} we obtain that $\tilde{\eta}^0_j(t)\pm \beta_{i,l} \tilde{\Phi}^i_{jl}(t)\geq 0$, where $\beta_{i,l}=\int_{T_{i-1}}^{T_i} \alpha_{i,l}(s) ds$. If $\beta_{i,l}=0$  then by the convention stated above $\tilde{\bPhi}^{i}_{\cdot,l}(t)=0$. Otherwise $\beta_{i,l}>0$
and thus $|\tilde{\Phi}^i_{jl}(t)|\leq 1/\beta_{i,l}$ for all $j\in [J]$. We can therefore conclude that $\tilde{\Phi}^i_{jl}$ is bounded in $L^\infty[0,T]$ and, since $T$ is finite, also in $L^2[0,T]$.
\end{proof}

Moreover, denoting by $B_{\W}$ and $B_{\U}$ uniform $L^2$-norm bounds on $\W$ and $\U$, respectively, and noting that $F$ is biaffine and that its linear and bilinear components define bounded linear operators and functionals with respect to the corresponding $L^2$-norms, there exists a constant $C>0$ such that
\begin{align}
\norm{\nabla_{\bw}F(\bw,\bepsilon)}_{L^2}
&\leq C\left(1+\norm{\bepsilon}_{L^2}\right)
\leq C(1+B_{\U})
=:G_{\W},
\label{eq:bound_grad_w}\\
\norm{\nabla_{\bepsilon}F(\bw,\bepsilon)}_{L^2}
&\leq C\left(1+\norm{\bw}_{L^2}\right)
\leq C(1+B_{\W})
=:G_{\U},
\label{eq:bound_grad_epsilon}
\end{align}
for all $(\bw,\bepsilon)\in\W\times\U$. Thus, $F(\cdot,\bepsilon)$ and $-F(\bw,\cdot)$ are uniformly Lipschitz continuous over $\W$ and $\U$, respectively, with Lipschitz constants $G_{\W}$ and $G_{\U}$.

For any $\bw\in\W$, define its robust objective value as
\begin{align}\label{eq:def_barJ}
\bar{J}(\bw)
=
\sup_{\bepsilon\in\U}F(\bw,\bepsilon),
\end{align}
and so
\(
\bar{J}^*
:=
\inf_{\bw\in\W}\bar{J}(\bw)
\)
denotes the optimal value of \eqref{eq:final_linear_RC}. For any
$(\bw,\bepsilon)\in\W\times\U$, define the saddle gap as
\[
\operatorname{Gap}(\bw,\bepsilon)
:=
\bar{J}(\bw)
-
\inf_{\bw'\in\W}F(\bw',\bepsilon).
\]
Since
\[
\inf_{\bw'\in\W}F(\bw',\bepsilon)
\leq
\bar{J}^*
\leq
\bar{J}(\bw),
\]
it follows that
\begin{equation}\label{eq:opt_gap_sad_gap}
0
\leq
\bar{J}(\bw)-\bar{J}^*
\leq
\operatorname{Gap}(\bw,\bepsilon).
\end{equation}
Consequently, if
$\operatorname{Gap}(\bw,\bepsilon)\leq\delta$,
then $\bw$ is a $\delta$-optimal solution of
\eqref{eq:final_linear_RC}.

We now introduce an SCLP-based variant of the SP--FW method, given in
Algorithm~\ref{alg:SCLP_SPFW}. Based on the online Frank--Wolfe algorithm of \cite{hazan2012projection}, the method generates feasible controller and uncertainty sequences. At each iteration, it requires the optimization of linear functionals over the controller and uncertainty feasible sets. We therefore define the corresponding linear optimization oracles, referred to as the $\W$-LMO and the $\U$-LMO, respectively.

For $\W$-LMO, let 
\[ \mathcal Q_{\W} = \left\{ \bb{q}\in\M^{[0,T],\dim(\bw)} : \bb{q}\text{ is piecewise affine on a finite partition of }[0,T] \right\}\]
be the coefficient class. For every $\bb{q}\in\mathcal Q_{\W}$ the corresponding oracle is defined as
\[
\operatorname{LMO}_{\W}(\bb{q})
\in
\argmin_{\bw\in\W}\int_{0}^T \bb{q}(t)^\top\bw(t)\,dt.
\]
Under the assumed availability of the SCLP--Simplex method, $\operatorname{LMO}_{\W}(\bb{q})$ can be computed by applying this method to the corresponding problem, yielding a piecewise-constant optimizer.

Moreover, for fixed $\bepsilon$, the coefficients of
$F(\cdot,\bepsilon)$ can be expressed in closed form whenever the
interval integrals of $\bepsilon$ can be computed. Thus, $\nabla_{\bw} F(\bw,\bepsilon)$ is affine in time for all $\bepsilon\in \U$ and therefore belongs to $\mathcal Q_{\W}$.

For $\U$-LMO, we restrict attention at this stage to the coefficient class
\[ \mathcal Q_{\U} := \left\{ \bb{q}\in\M^{[0,T],K} : \sup_{t\in[0,T]}\norm{\bb{q}(t)}<\infty \right\}.\]
We assume that, for every $\bb{q}\in\mathcal Q_{\U}$,  the oracle 
\[ \operatorname{LMO}_{\U}(\bb{q}) \in \argmax_{\bepsilon\in\U} \int_0^T \bb{q}(t)^\top\bepsilon(t)\,dt, \]
is well-defined. In Sections \ref{sec:U1_LMO} and \ref{sec:5.2}, we establish the existence and computability of this oracle for the exact formulation under $\U_1$ and the approximate formulation under $\U_3$.

\begin{algorithm}[t]
\caption{SP--FW algorithm}
\label{alg:SCLP_SPFW}
\begin{algorithmic}[1]
\State {\bf Input:} $\sigma_{\W}^{(1)},\sigma_{\U}^{(1)}\in \Real_+$, piecewise-constant $\bw^{(1)}\in\W$ and
$\bepsilon^{(1)}\in\U$, and set
$\bz^{(1)}=\bar{\bz}^{(1)}
:=(\bw^{(1)},\bepsilon^{(1)})$. 
\For{$\ell=1,2,\ldots$}
\State Set
    \vspace{-10pt}\[
    {\bb q}_{\W}^{(\ell)}
    :=
    \frac{1}{\ell}
    \sum_{r=1}^{\ell}
    \nabla_{\bw}
    F\!\left(\bw^{(r)},\bepsilon^{(r)}\right)
    +
    \sigma^{(\ell)}_{\W}
    \left(\bw^{(\ell)}-\bw^{(1)}\right),
    \]
    and
    \vspace{-10pt}\[
    {\bb q}_{\U}^{(\ell)}
    :=
    \frac{1}{\ell}
    \sum_{r=1}^{\ell}
    \nabla_{\bepsilon}
    F\!\left(\bw^{(r)},\bepsilon^{(r)}\right)
    -
    \sigma^{(\ell)}_{\U}
    \left(\bepsilon^{(\ell)}-\bepsilon^{(1)}\right).
    \]
    \State Compute the Frank--Wolfe directions \vspace{-10pt}
    \begin{align*}
    \tilde{\bw}^{(\ell)}
    =
    \operatorname{LMO}_{\W}
    \!\left({\bb q}_{\W}^{(\ell)}\right),
    \qquad
    \tilde{\bepsilon}^{(\ell)}
    =
    \operatorname{LMO}_{\U}
    \!\left({\bb q}_{\U}^{(\ell)}\right).
    \end{align*}
    \vspace{-15pt}
    \State Set $\gamma_{\ell}:=\ell^{-1/4}$ and
    $\bs^{(\ell)}
    :=(\tilde{\bw}^{(\ell)},\tilde{\bepsilon}^{(\ell)})$.
    \State Update $\sigma^{(\ell+1)}_{\W}=\frac{\sigma^{(\ell)}_{\W}\ell+\sigma_{\W}^{(1)}(\ell+1)^{-1/4}}{\ell+1},\; \sigma^{(\ell+1)}_{\U}=\frac{\sigma^{(\ell)}_{\U}\ell+\sigma_{\U}^{(1)}(\ell+1)^{-1/4}}{\ell+1}$
    \State Update
    \(
    \bz^{(\ell+1)}
    :=
    (1-\gamma_{\ell})\bz^{(\ell)}
    +
    \gamma_{\ell}\bs^{(\ell)}.
    \)
    \State Update
    \(
    \bar{\bz}^{(\ell+1)}
    :=
    \frac{\ell}{\ell+1}\bar{\bz}^{(\ell)}
    +
    \frac{1}{\ell+1}\bz^{(\ell+1)}.
    \)
\EndFor
\end{algorithmic}
\end{algorithm}

Given these oracles the algorithm is well defined and the following result establishes, by extending the regret analysis of \cite{hazan2012projection} to the present setting, that the ergodic averages of these sequences have a vanishing gap, which provides an optimality certificate for the averaged controller.
\begin{theorem}\label{thm:SP_FW_convergence}
Let $\sigma_{\W}^{(1)}=2G_{\W}/D_{\W}$ and
$\sigma_{\U}^{(1)}=2G_{\U}/D_{\U}$ be the regularization inputs to Algorithm~\ref{alg:SCLP_SPFW}, 
where $D_{\W}$ and $D_{\U}$ are the diameters of $\W$ and $\U$,
respectively, under the corresponding $L^2$-norms, and $G_{\W}$ and
$G_{\U}$ are defined in
\eqref{eq:bound_grad_w}--\eqref{eq:bound_grad_epsilon}.
Then, Algorithm~\ref{alg:SCLP_SPFW} is well-defined and generates a sequence $\{\bz^{(\ell)}\}_{\ell\in\mathbb{N}}$, where $\bz^{(\ell)}=(\bw^{(\ell)},\bepsilon^{(\ell)})$. 

For any $L\in \mathbb{N}$, let $\bar{\bz}^{(L)} = \bigl(\bar{\bw}^{(L)},\bar{\bepsilon}^{(L)}\bigr)$ be the corresponding ergodic average generated by the algorithm at iteration $L$. Then, $\bar{\bw}^{(L)}\in\W,$ is piecewise constant, $
\bar{\bepsilon}^{(L)}\in\U$, and 
\begin{align}
\bar{J}(\bar{\bw}^{(L)})-\bar{J}^{*}\leq \operatorname{Gap}
\bigl(\bar{\bw}^{(L)},\bar{\bepsilon}^{(L)}\bigr)
\leq
57\left(
D_{\W}G_{\W}
+
D_{\U}G_{\U}
\right)
L^{-1/4}.
\label{eq:SP_FW_gap}
\end{align}
\end{theorem}

\begin{proof}
By assumption, $\bw^{(1)}\in\W$ is piecewise constant and $\bepsilon^{(1)}\in\U$. Suppose that $\bw^{(r)}$ is piecewise constant on a finite partition of $[0,T]$ for every $r\in[\ell]$. By \eqref{eqn:objective_first_form}, for every $\bepsilon\in\U$, $\nabla_{\bw}F(\bw,\bepsilon)$ is affine in time. Moreover, $\bw^{(\ell)}-\bw^{(1)}$ is piecewise constant on a finite partition. Hence, ${\bb q}_{\W}^{(\ell)}$ is piecewise affine on a finite common refinement, and therefore ${\bb q}_{\W}^{(\ell)}\in\mathcal Q_{\W}.$ Similarly, since $F$ is affine in $\bepsilon$ and $\W$ is uniformly bounded, $\nabla_{\bepsilon}F(\bw^{(r)},\bepsilon^{(r)})$ is bounded for every $r\in[\ell]$. Together with the uniform boundedness of $\U$, this implies that ${\bb q}_{\U}^{(\ell)}\in\mathcal Q_{\U}$. Thus, both oracle calls at iteration $\ell$ are well-defined. 

Since the oracles return elements of $\W$ and $\U$, respectively, and both sets are convex, the Frank--Wolfe updates imply that $\bw^{(\ell+1)}\in\W$ and $\bepsilon^{(\ell+1)}\in\U$. In addition, $\bw^{(\ell+1)}$ is piecewise constant on a finite common refinement. This completes the induction. Furthermore, the ergodic averages are convex combinations of feasible iterates. Thus, $\bar{\bw}^{(L)}\in\W,\;
\bar{\bepsilon}^{(L)}\in\U$ for every $L\in\mathbb N$.

    Using the construction in the proof of~\cite[Theorem~4.4]{hazan2012projection}, Algorithm~\ref{alg:SCLP_SPFW} applies the online FW algorithm with exponent $a=1/4$ to the  regularized surrogate losses associated with functions
    $f_\ell(\bw)=F(\bw,{\bepsilon}^{(\ell)})$
    and $g_{\ell}(\bepsilon)=-F({\bw}^{(\ell)},\bepsilon)$. By \eqref{eq:bound_grad_w}--\eqref{eq:bound_grad_epsilon}, these functions are Lipschitz continuous for every $\ell$, with Lipschitz constants $G_{\W}$ and $G_{\U}$, respectively.
    The proof of the regret guarantee in \cite[Theorem~4.4]{hazan2012projection} extends to the present infinite-dimensional setting by carrying out the same argument in the $L^2$-norm, with respect to the corresponding infima. The argument requires that the underlying feasible set be convex and bounded, that the online loss functions be convex and Lipschitz continuous with respect to the $L^2$-norm, and that the required linear minimization oracle be well-defined. Since these conditions hold in our setting, the regret bounds of \cite[Theorem~4.4]{hazan2012projection} apply to the original loss functions.
    
    Using the biaffine nature of $F$ we obtain that
    \begin{align*}
        &\frac{1}{L}\left(\sum_{\ell=1}^{L} f_\ell(\bw^{(\ell)})-\inf_{\bw\in\W} \sum_{\ell=1}^{L} f_\ell(\bw)\right)%
       = \frac{1}{L}\sum_{\ell=1}^{L} F(\bw^{(\ell)},\bepsilon^{(\ell)})-\inf_{\bw\in\W} F(\bw,\bar{\bepsilon}^{(L)}),\\
        &\frac{1}{L}\left(\sum_{\ell=1}^{L} g_\ell(\bepsilon^{(\ell)})-\inf_{\bepsilon\in\U} \sum_{\ell=1}^{L} g_\ell(\bepsilon)\right)= \sup_{\bepsilon\in\U} F(\bar{\bw}^{(L)},\bepsilon)-\frac{1}{L}\sum_{\ell=1}^{L} F(\bw^{(\ell)},\bepsilon^{(\ell)}).
    \end{align*}
    Summing the above equalities and applying the regret bounds together with \eqref{eq:opt_gap_sad_gap}   
    results in inequality \eqref{eq:SP_FW_gap}. 
\end{proof}

Theorem~\ref{thm:SP_FW_convergence} establishes that the optimality gap of the averaged controller generated by Algorithm~\ref{alg:SCLP_SPFW} converges to zero. However, as stated above, computing the $\operatorname{LMO}_{\U}\bigl({\bb q}_{\U}^{(\ell)}\bigr)$ at every iteration may be computationally challenging under the general uncertainty set $\U_3$. We therefore consider two settings. In Section~\ref{sec:U1_LMO}, we show that the oracle can be computed exactly under $\U_1$. In Section~\ref{sec:5.2}, we introduce an approximation of $F$ that yields a tractable uncertainty-side oracle under $\U_3$.

\subsection{Computing the $\U$-LMO under $\U_1$} \label{sec:U1_LMO}

In this section, we consider the case of budgeted uncertainty $\U_1$. We start by showing that computing the $\U$-LMO can be done efficiently.

Let $\bb q:[0,T]\rightarrow\mathbb R^K$ be piecewise affine on a finite partition of $[0,T]$. Using Definition~\ref{def:avop}, we define the permutation-valued function $\tilde{\bpi}(t):=\pi(\bb q(t))$, where $\tilde{\bpi}(t)$ is the absolute-value ordering permutation of $\bb q(t)$. Since $\bb q$ is piecewise affine on a finite partition and $\Gamma$ is piecewise constant by Assumption~\ref{ass:Gamma}, this partition can be refined to obtain a finite ordered partition $(\tilde{T}_i)_{i=0}^{\tilde{p}}$ such that $\tilde{\bpi}(\cdot)$, $\sign(\bb q(\cdot))$, and $\Gamma(\cdot)$ are constant on the interior of each partition interval. 
Thus, the $\U$-LMO can be computed pointwise on each interval, as shown in the following lemma.

\begin{lemma}\label{lem:U1_LMO} Let $\bb q:[0,T]\rightarrow\mathbb R^K$ be piecewise affine on a finite partition, and let $(\tilde{T}_i)_{i=0}^{\tilde{p}}$ be the ordered partition of $[0,T]$ described above. For each $i\in[\tilde{p}]$, let $\tilde{\bpi}^{i}$,  $\tilde{\bs}^{i}$, and $\hat{\Gamma}_i$ denote the constant values of $\bpi(\bb q(\cdot))$, ${\sign}(\bb q(\cdot))$, and $\lfloor\Gamma(\cdot)\rfloor$, respectively, on the interior of $(\tilde{T}_{i-1},\tilde{T}_i]$. Then, the piecewise-constant function $\bepsilon^\star$ defined for all $i\in[\tilde{p}]$ and $t\in(\tilde{T}_{i-1},\tilde{T}_i)$ by
\[ \epsilon_k^\star(t) = \begin{cases} \tilde{s}^{i}_k, & k\in \left\{ \tilde{\pi}^{i}_1,\ldots, \tilde{\pi}^{i}_{\hat{\Gamma}_i} \right\}, \\[0.3em] \tilde{s}^{i}_k \bigl(\Gamma(t)-\hat{\Gamma}_i\bigr), & \hat{\Gamma}_i<K \text{ and } k=\tilde{\pi}^{i}_{\hat{\Gamma}_i+1}, \\[0.3em] 0, & \text{otherwise}, \end{cases} \]
with arbitrary feasible values at the partition points, satisfies $\bepsilon^\star \in \operatorname{LMO}_{\U}(\bb q)$. \end{lemma}
The result follows directly from the construction of the partition and the standard solution of the budgeted uncertainty problem in~\cite{Bertsimas2004}.

We next verify that the preceding lemma applies to the $\U$-LMO at every iteration of Algorithm~\ref{alg:SCLP_SPFW}. Suppose that $\bepsilon^{(1)}$ is piecewise constant. Theorem~\ref{thm:SP_FW_convergence} establishes that each controller iterate $\bw^{(r)}$ is piecewise constant on a finite partition. Hence, by~\eqref{eq:objective_first_form}, $\nabla_{\bepsilon}F(\bw^{(r)},\bepsilon^{(r)})$ is piecewise affine for every $r$. If $\bepsilon^{(\ell)}$ is piecewise constant, then the regularization term $\bepsilon^{(\ell)}-\bepsilon^{(1)}$ is also piecewise constant. Therefore, $\bb q_{\U}^{(\ell)}$ is piecewise affine on a finite common refinement. The preceding lemma then yields a piecewise-constant solution $\tilde{\bepsilon}^{(\ell)}\in\operatorname{LMO}_{\U}(\bb q_{\U}^{(\ell)})$. Since the Frank--Wolfe update is a convex combination of $\bepsilon^{(\ell)}$ and $\tilde{\bepsilon}^{(\ell)}$, the next iterate $\bepsilon^{(\ell+1)}$ is also piecewise constant on a finite common refinement. Thus, by induction, the $\U$-LMO can be computed using the preceding lemma at every iteration of Algorithm~\ref{alg:SCLP_SPFW}.

Moreover, since the resulting functions \(\bepsilon\) are piecewise constant, the coefficients of \(\bw\) in the \(\W\)-LMO can be computed efficiently by replacing the integrals in~\eqref{eqn:objective_first_form} with finite sums over constant values.

\subsection{Objective Approximation}\label{sec:5.2}
For the general uncertainty set $\U_3$, it is not clear whether the $\U$-LMO can be solved tractably. We therefore consider applying Algorithm~\ref{alg:SCLP_SPFW} to an approximation of $F$. 
For each $p\in\mathbb{N}$, consider the uniform partition of $[0,T]$ given by $T_i=iT/p$, $i\in[p]_0$, and let $\W_p$ denote the feasible controller set of~\eqref{eq:final_linear_RC} associated with this partition.
In this subsection, $F$ and $\tilde{F}$ defined in~\eqref{eqn:objective_first_form} are understood as functionals on $\W_p\times\U_3$.

Since any $\bepsilon\in\U_3$ satisfies the constraint $ |\bepsilon(s)| \le 1$ for all $s\in[0,T]$, we obtain %
\begin{equation}\label{eq:tail_to_zero}
\Big| (\tilde{\blambda} \circ \bc)^\top \left(\sum_{i=1}^p \int_{T_{i-1}}^{T_i} \int_{T_{i-1}}^t \bepsilon(s) \, ds \, dt\right)\Big|
\;\le\; \frac{T^2 \norm{\tilde{\blambda} \circ \bc}_{1}}{2p}
\;\xrightarrow[p\to\infty]{} 0.\end{equation}
This result motivates omitting the rightmost term in~\eqref{eqn:objective_first_form} and replacing $F$ with the approximation $\tilde{F}$.

Using $\tilde{F}$ in Algorithm~\ref{alg:SCLP_SPFW}, the corresponding $\U$-LMO can be solved by applying Theorem~\ref{thm:robust_counterpart_constraint_in_3} and complementarity to recover variables $\bxi_{i,n}$ satisfying
$S_{i,n}\bxi_{i,n}
=
\int_{R_{i,n}}\bepsilon(s)\,ds$ and
$\sum_{n=0}^{K}S_{i,n}\bxi_{i,n}
=
\int_{T_{i-1}}^{T_i}\bepsilon(s)\,ds$.
Let $\Xi$ denote the set of all vectors $\bxi=(\bxi_{i,n})_{i\in[p],\,n\in[K]_0}$ obtained from some $\bepsilon\in\U_3$ under this mapping. Since $\tilde{F}$ depends on $\bepsilon$ only through these aggregate quantities, the objective
$\hat{F}(\bw,\bxi):=\tilde{F}(\bw,\bepsilon)$
is well defined. Thus, Algorithm~\ref{alg:SCLP_SPFW} can be applied directly with $\bxi\in\Xi$ in place of $\bepsilon\in\U_3$.

While it is clear that the values of $F$ and $\tilde{F}$ are close for a fixed $\bw$, the following theorem establishes that this approximation gap also carries over to the optimal values obtained using the two objective functions.
\begin{theorem}[Consistency] \label{thm:F2_approx} For any $p\in \mathbb{N}$ let \(f_p^{1}, f_p^{2} : \W_p \times \U \to \Real\) be functionals satisfying
\begin{equation}\label{eq:cond}
\left|f_p^{1}(\bw,\bepsilon)-f_p^{2}(\bw,\bepsilon)\right| \le \frac{C}{p}
\qquad \forall\, \bw \in \W_p,\ \bepsilon \in \U,
\end{equation}
for some constant \(C>0\) independent of \(p\). For \(i \in \{1,2\}\), define the max-envelopes
\[
\mathcal{F}_p^{\,i}(\bw) \triangleq \sup_{\bepsilon \in \U} f_p^{\,i}(\bw,\bepsilon),
\qquad \bw \in \W_p,
\]
and assume $\mathcal{F}_p^{\,i}(\bw)$ is finite and uniformly bounded from below for all \(\bw \in \W_p\) and all \(p \in \mathbb{N}\). 
Then the associated sup-envelopes satisfy
\begin{equation}\label{eq:diff_abs}
\bigl| \mathcal{F}_{p}^1(\bw)
- \mathcal{F}_{p}^2(\bw) \bigr|
\leq \frac{C}{p}, \quad \forall \bw \in \W_p.
\end{equation}
Consequently, their infimal values are also \(O(1/p)\)-close, namely,
\begin{equation}\label{eq:diff_opt_abs}
\left|\inf_{\bw\in \W_p}\mathcal{F}_p^{\,1}(\bw)-\inf_{\bw\in \W_p}\mathcal{F}_p^{\,2}(\bw)\right|
\le \frac{C}{p}.
\end{equation}

\end{theorem}

\begin{proof}
Fix \(p\in\mathbb N\). By \eqref{eq:cond}, for every \(\bw\in\W_p\) and \(\bepsilon\in\U\),
\[
f_p^2(\bw,\bepsilon)-\frac{C}{p} \leq f_p^1(\bw,\bepsilon)\le f_p^2(\bw,\bepsilon)+ \frac{C}{p}.
\]
Taking the supremum over \(\bepsilon\in\mathcal U\) gives
\begin{equation}\label{eq:diff_non_abs}
\mathcal{F}_{p}^2(\bw) - \frac{C}{p} 
\leq \mathcal{F}_{p}^1(\bw) \leq   \mathcal{F}_{p}^2(\bw) + \frac{C}{p}.
\end{equation}
Thus, \eqref{eq:diff_abs} holds.
Taking the infimum over $\bw\in \W_p$ in \eqref{eq:diff_non_abs},  establishes \eqref{eq:diff_opt_abs}. %
\end{proof}

In our setting, $f^{1}_p(\cdot)=F(\cdot)$ and $f^{2}_p(\cdot)=\tilde{F}(\cdot)$, resulting in $\mathcal{F}^{1}_p(\cdot)=\bar{J}_p(\cdot)$ given by \eqref{eq:def_barJ}, and $\mathcal{F}_p^{2}(\cdot)=\tilde{J}_p(\cdot)$, defined analogously. %
We denote the corresponding infima over $\W_p$ by $\bar{J}^*_p$ and $\tilde{J}^*_p$, respectively. Utilizing \eqref{eq:tail_to_zero} together with Theorem~\ref{thm:F2_approx} we derive the following guarantee on the performance of a controller obtained from the approximation $\tilde{F}$.
\begin{corollary}[Performance of Near-Optimal Solution]
Let $C=\frac{ T^2 \norm{\tilde{\blambda} \circ \bc}_1}{2}$.
Then, for any $p\in\mathbb{N}$ and any $\delta>0$, any %
$\tilde{\bw}^*_{p,\delta}\in \W_p$ such that $\tilde{J}_p(\tilde{\bw}^*_{p,\delta})-\tilde{J}^*_p\leq \delta$ satisfies %
\[
\bar{J}(\tilde{\bw}^*_{p,\delta}) - \bar{J}^*_p  \le \frac{2C}{p}+\delta.
\]
In particular, for any prescribed tolerance $\varepsilon>0$, choosing $p$ sufficiently large and $\delta>0$ sufficiently small so that
$\frac{2C}{p}+\delta\leq\varepsilon$
ensures that
$\bar{J}(\tilde{\bw}_{p,\delta}^*)-\bar{J}^*_p\leq\varepsilon$.
\end{corollary}

\section{Conclusion}
In this work, we introduced a partially adaptive robust optimization framework for uncertain fluid processing networks, with the aim of informing control decisions in their discrete counterparts. The proposed policies adapt to network states observed at predefined times. We considered three types of uncertainty sets that capture bounds on adversarial noise over time and across buffers, and developed SCLP formulations that yield either exact solutions or highly accurate approximations through tailored solution procedures based on the SCLP--Simplex algorithm.

Future work will focus on extending the SCLP--Simplex algorithm to accommodate the additional constraints imposed by these formulations. This extension would enable efficient implementation of the proposed framework and support the derivation of robust and adaptive control policies for discrete processing networks.

\bibliographystyle{plainnat}
\bibliography{bibliography}%

\begin{thebibliography}{15}
\providecommand{\natexlab}[1]{#1}
\providecommand{\url}[1]{\texttt{#1}}
\expandafter\ifx\csname urlstyle\endcsname\relax
  \providecommand{\doi}[1]{doi: #1}\else
  \providecommand{\doi}{doi: \begingroup \urlstyle{rm}\Url}\fi

\bibitem[Anderson et~al.(1983)Anderson, Nash, and
  Perold]{AndersonNashPerold1983}
Edward~J. Anderson, Peter Nash, and Andr{\'e}~F. Perold.
\newblock Some properties of a class of continuous linear programs.
\newblock \emph{SIAM J. Control Optim.}, 21\penalty0 (5):\penalty0 758--765,
  1983.
\newblock \doi{10.1137/0321046}.

\bibitem[Ben-Tal et~al.(2004)Ben-Tal, Goryashko, Guslitzer, and
  Nemirovski]{BenTal2004}
A.~Ben-Tal, A.~Goryashko, E.~Guslitzer, and A.~Nemirovski.
\newblock Adjustable robust solutions of uncertain linear programs.
\newblock \emph{Math. Program.}, 99:\penalty0 351--376, 2004.
\newblock \doi{10.1007/s10107-003-0454-y}.

\bibitem[Ben-Tal et~al.(2009)Ben-Tal, {El Ghaoui}, and Nemirovski]{BenTal2009}
Aharon Ben-Tal, Laurent {El Ghaoui}, and Arkadi Nemirovski.
\newblock \emph{Robust Optimization}.
\newblock Princeton Series in Applied Mathematics. Princeton University Press,
  Princeton, NJ, 2009.
\newblock ISBN 978-1-4008-3105-0.
\newblock \doi{10.1515/9781400831050}.

\bibitem[Bertsimas and Sim(2004)]{Bertsimas2004}
Dimitris Bertsimas and Melvyn Sim.
\newblock The price of robustness.
\newblock \emph{Oper. Res.}, 52\penalty0 (1):\penalty0 35--53, 2004.
\newblock \doi{10.1287/opre.1030.0065}.

\bibitem[Bertsimas et~al.(2015)Bertsimas, Nasrabadi, and
  Paschalidis]{Bertsimas2015}
Dimitris Bertsimas, Ebrahim Nasrabadi, and Ioannis~Ch. Paschalidis.
\newblock Robust fluid processing networks.
\newblock \emph{IEEE Trans. Autom. Control}, 60\penalty0 (3):\penalty0
  715--728, 2015.
\newblock \doi{10.1109/tac.2014.2352711}.

\bibitem[Ghate(2020)]{Ghate2020}
Archis Ghate.
\newblock Robust continuous linear programs.
\newblock \emph{Optim. Lett.}, 14\penalty0 (5):\penalty0 1627--1642, 2020.
\newblock \doi{10.1007/s11590-020-01539-6}.

\bibitem[Hazan and Kale(2012)]{hazan2012projection}
Elad Hazan and Satyen Kale.
\newblock Projection-free online learning.
\newblock In \emph{Proceedings of the 29th International Conference on Machine
  Learning (ICML)}, pages 1843--1850, 2012.

\bibitem[Nazarathy and Weiss(2010)]{Nazarathy2009}
Yoni Nazarathy and Gideon Weiss.
\newblock Positive {H}arris recurrence and diffusion scale analysis of a push
  pull queueing network.
\newblock \emph{Perform. Eval.}, 67\penalty0 (4):\penalty0 201--217, 2010.
\newblock \doi{10.1016/j.peva.2009.09.010}.

\bibitem[Pullan(1993)]{Pullan1993}
Malcolm~C. Pullan.
\newblock An algorithm for a class of continuous linear programs.
\newblock \emph{SIAM J. Control Optim.}, 31\penalty0 (6):\penalty0 1558--1577,
  1993.
\newblock \doi{10.1137/0331073}.

\bibitem[Shapiro(2001)]{Shapiro2001}
Alexander Shapiro.
\newblock On duality theory of conic linear problems.
\newblock In M.~A. Goberna and M.~A. L{\'o}pez, editors, \emph{Semi-Infinite
  Programming: Recent Advances}, volume~57 of \emph{Nonconvex Optimization and
  Its Applications}, pages 135--165. Kluwer Academic Publishers, Dordrecht,
  2001.
\newblock \doi{10.1007/978-1-4757-3403-4_7}.

\bibitem[Shindin et~al.(2021)Shindin, Masin, Weiss, and
  Zadorojniy]{Shindin2021}
Evgeny Shindin, Michael Masin, Gideon Weiss, and Alexander Zadorojniy.
\newblock Revised {SCLP}-simplex algorithm with application to large-scale
  fluid processing networks.
\newblock In \emph{2021 IEEE 60th Conference on Decision and Control (CDC)},
  pages 3863--3868, Piscataway, NJ, USA, 2021. IEEE.
\newblock \doi{10.1109/CDC45484.2021.9683390}.

\bibitem[Ship et~al.(2022)Ship, Shindin, Boni, and Dattner]{Ship2022}
Harold Ship, Evgeny Shindin, Odellia Boni, and Itai Dattner.
\newblock A robust server-effort policy for fluid processing networks.
\newblock In \emph{2022 IEEE 61st Conference on Decision and Control (CDC)},
  pages 5902--5909, Piscataway, NJ, 2022. IEEE.
\newblock \doi{10.1109/CDC51059.2022.9992816}.

\bibitem[Weiss(1990)]{Weiss1990}
Gideon Weiss.
\newblock A class of continuous linear programming problems with state
  constraints.
\newblock \emph{SIAM J. Control Optim.}, 28\penalty0 (2):\penalty0 312--323,
  1990.
\newblock \doi{10.1137/0328017}.

\bibitem[Weiss(2008)]{Weiss2008}
Gideon Weiss.
\newblock A simplex based algorithm to solve separated continuous linear
  programs.
\newblock \emph{Math. Program.}, 115\penalty0 (1):\penalty0 151--198, 2008.
\newblock \doi{10.1007/s10107-008-0217-x}.

\bibitem[Z{\u{a}}linescu(2002)]{zalinescu2002convex}
Constantin Z{\u{a}}linescu.
\newblock \emph{Convex Analysis in General Vector Spaces}.
\newblock World Scientific, River Edge, NJ, 2002.
\newblock ISBN 978-981-238-067-8.
\newblock \doi{10.1142/5021}.

\end{thebibliography}

\end{document}